\documentclass[11pt,english,a4paper]{amsart}

\usepackage{babel}
\usepackage{amsfonts,amssymb,latexsym,xspace,epsfig,graphics,color}  
\usepackage{amsmath,enumerate}  
\numberwithin{equation}{section}

\title[Exponential inequalities and functional estimation]{Exponential inequalities and functional estimations for weak dependent datas~; applications to dynamical systems.}
\author{ V. Maume-Deschamps}
\address{Universit\'e de Bourgogne   
              B.P. 47870   
         21078 Dijon Cedex FRANCE}
\email{vmaume@u-bourgogne.fr}   

\thanks{ \ }

 \newcommand{\eps}{\varepsilon}   
\newcommand{\N}{\ensuremath{\mathbb{N}}}    
\newcommand{\R}{\ensuremath{\mathbb{R}}}     
\newcommand{\Z}{\ensuremath{\mathbb{Z}}}    
\renewcommand{\P}{\ensuremath{\mathbb{P}}}       
\newcommand{\E}{\ensuremath{\mathbb{E}}}

\newcommand{\Id}{\mathbf {1}}

\newcommand{\MM}{\mathcal {M}}

\newcommand{\ZZ}{\mathcal Z}
 
\newcommand{\CC}{\mathcal C}

 \newcommand{\PP}{\mathcal P}

\newcommand{\Lip}{\mbox{\rm Lip}}

\renewcommand{\hat}{\widehat}

\newtheorem{theo}{Theorem}[section]
\newtheorem*{theointro}{Theorem}  
\newtheorem*{corointro}{Corollary}  
\newtheorem{prop}[theo]{Proposition}    
\newtheorem{coro}[theo]{Corollary}    
\newtheorem{lem}[theo]{Lemma}

\theoremstyle{definition}
\newtheorem{defi}{Definition}
 
\newtheorem{cond}{Assumption}
\theoremstyle{remark}
\newtheorem*{rem}{Remark}
\newtheorem*{rems}{Remarks}

\begin{document}
\begin{abstract}
We estimate density and regression functions for weak dependant datas.  Using an exponential inequality obtained in \cite{DP} and  in \cite{taudep}, we control the deviation between the estimator and the function itself. These results are applied to a large class of dynamical systems and lead to estimations of invariant densities and on the mapping itself.
\end{abstract}
\maketitle
Dynamical systems are widely used by scientists to modalize complex systems (\cite{datas}). Therefore, estimating functions related to dynamical systems is crucial. Of particular interest are~: the invariant density, the mapping itself, the pressure function. We shall see that  many dynamical systems have the same behavior as weak dependant processes (as defined in \cite{DL}).  We obtain results of deviation for regression functions and densities for weak dependant processes and apply these results to dynamical systems. \\
In \cite{taudep} we gave an estimation (with control of the deviation) of the pressure function for some expanding maps of the interval. Results on the estimation of the invariant density for the same kind of maps where obtained in \cite{clem} and stated in \cite{maes}. In this later article, results on the estimation of the mapping were also stated. These last two papers dealt with convergence in quadratic mean. Our goal here is, on one hand, to consider more general dynamical systems~: non uniformly hyperbolic maps on the interval, dynamics in higher dimension ... On the second hand, we obtain bounds on the deviation between the estimator and the regression function as well as almost sure convergence. Related results on the estimation of the regression function may also be found in \cite{FV} and \cite{Mas} where strongly mixing processes are considered and almost everywhere convergence and asymptotic normality are proved. Our aim is to provide an estimation of the deviation between the estimator and the regression function for a larger class of mixing processes and for regression functions that may have singularities. \\
Before giving the precise definitions and results, let us state our main results informally. 
We consider a weak dependant  stationary process  $X_0$, ..., $X_i$, ...  taking values in $\Sigma \subset \R^d$. Our condition of weak dependence is  with respect to a Banach space $\CC$ of bounded functions on $\Sigma$ (see Definition \ref{weakmixing} below). Let $(Y_i)_{i\in\N}$ be a  stationary process taking values in $\R$ and satisfying a condition of weak dependence according with the process  $X_0$, ..., $X_i$, ... Consider the regression function $r(x)=\E(Y_i\ | \ X_i=x)$, we shall assume some regularity on $r$ (see Assumption \ref{almosthol} below). We shall also consider that the process $(X_i)_{i\in\N}$ has a density $f$ with respect to the Lebesgue measure $m$. $f$ is assmued to have  some regularity properties that allow localised singularities. \\
Consider a non negative kernel $K\in\CC$ satisfying Assumption \ref{Kernel}  and the estimators (introduced in \cite{Ro} and used for example in  \cite{clem}, \cite{maes}, \cite{FV}, \cite{Mas})~:
\begin{eqnarray}
\widehat{f}_n(x) =  \frac1{nh} \sum_{i=0}^{n-1} K\left(\frac{x-X_i}h\right) && \widehat{g}_n(x) = \frac1{nh} \sum_{i=0}^{n-1} Y_i \/ K\left(\frac{x-X_i}h\right) \nonumber\\
\lefteqn{ \widehat{r}_n(x) = \frac{\widehat{g}_n(x)}{\widehat{f}_n(x)} \/. \label{defreg}}
\end{eqnarray}
Remark that $\hat{f}_n(x)=0$ implies $\hat{g}_n(x)=0$, in that case, we define $\hat{r}_n(x)=0$.
\begin{theointro}
Assume some weak dependence condition on $X_i$ and $Y_i$ (see Definition \ref{weakmixing}), assume  that  $\inf f>0$. There exists $M>0$, $L>0$, $R>0$, $0\leq \beta <1$, $D>0$, $\gamma'>0$ such that,  outside a set of measure less than $D h^{\gamma'}$, for all $n\in \N^\star$,  for all $t\in\R^{+}$,   for all $u\geq\mbox{Ct}\/ h$, 
$$\P(|\widehat{f}_n(x)-f(x)|>t-u^\alpha)\leq 2e^{\frac1e} \exp[-t^2 M h^{\beta+2} n] \ \mbox{if} \ f \ \mbox{is} \ \alpha \ \mbox{regular}$$
if $r$ is bounded and $\alpha$-regular,
$$\P(|\widehat{r}_n(x) -r(x)|>t-u^\alpha)\leq e^{\frac1e} \left(2\exp[-t^2 L n h^{\beta+2}]+ \exp[-L'\/n\/ h^{\beta+2}]\right) \/,$$
if $Y_i\in L^\infty$ and $r$ is $\alpha$-regular,
$$\P(|\widehat{r}_n(x) -r(x)|>t-u^\alpha)\leq 2e^{\frac1e} \exp[-t^2 L n h^{\beta+2}] \/.$$
As a consequence,  provided $h=h_n$ goes to zero and $n h_n^{\beta+2} =O(n^\eps)$, $\eps>0$, we obtain the following convergences provided $r$ (and $f$) are $\alpha$-regular~:
\begin{itemize}
\item for $m$ almost all $x\in \Sigma$, $\hat{f}_n(x)$ converges to $f(x)$ and $\hat{r}_n(x)$ converges to $r(x)$ almost surely and in $L^p$ for any $1\leq p$,
\item  $\displaystyle\E\left(\int_\Sigma |\hat{f}_n(x) -f(x)| \/dx \right)$ and $\displaystyle\E\left(\int_\Sigma |\hat{r}_n(x)-r(x)|\right)$ go to zero.
\item for almost all $x\in\Sigma$, for $a <\frac12$, $|\hat{g}_n(x)-g(x)| =O_\P(n^{-a})$ where $\hat{g}_n$ is either $\hat{f}_n$ or $\hat{r}_n$ and $g$ is either $f$ or $r$.
\end{itemize}
\end{theointro}
Using a ``time inverted process'' (as in \cite{DP,taudep}), we deduce an estimation of the invariant density and of the mapping itself for dynamical systems. Consider a discrete dynamical system from a ``good class'' (see Definition \ref{dyn} for precise settings), denote by $T$ the mapping and $f$ the invariant density of interest. Pick $X_0$ at random with stationary law $\mu=f \/ m$ and let $X_i=T^i(X_0)$. Consider a kernel $K$ and $\widehat{f}_n$ as above, let $\hat{r}_n^j$ be the estimator for   $Y_i^j=X_{i+1}^j+\eps_i^j$, $j=1\/, \ldots \/, d$, $X_i^j$ is the $j$th coordinate of $X_i$ and $(\eps_i^j)_{i\in\N}$ are independent stationary process with $0$ mean, independent of $(X_i)_{i\in\N}$. Let $\widehat{T}_n$ be an estimator of the map $T$ which has coordinates $\hat{r}_n^j$.  Let $\Vert \ \Vert_d$ be the sup norm in $\R^d$, we have the following results. 
\begin{corointro}
Assume that the dynamical is in the good class. There exists $M>0$, $L>0$, $R>0$, $0\leq \beta \leq 1$, $\gamma'>0$  such that outside a set of measure less than $R h^{\gamma'}$,   for all $t\in\R^{+}$,   for all $u\geq \mbox{Ct} \/h$, 
\begin{align}
\mbox{is $f$ is regular}\ \P(|\hat{f}_n(x)-f(x)|>t-u^\alpha)\leq 2e^{\frac1e} \exp[-t^2 M h^{\beta+2} n] \\
\P(\Vert\hat{T}_n(x)-T(x)\Vert_d>t-u^\alpha) \leq 2e^{\frac1e} \exp[-t^2 L h^{\beta+2}  n] \/.
\end{align}
As a consequence,  provided $h=h_n$ goes to zero  and $n h_n^{\beta+2} =O(n^\eps)$, $\eps>0$, we obtain the following convergences~:
\begin{itemize}
\item for $m$ almost all $x\in \Sigma$,  $\hat{T}_n(x)$ converges to $T(x)$ almost surely and in $L^p$  for any $1\leq p$, the same holds for $\hat{f}_n(x)$ provided that $f$ is $\alpha$-regular.
\item $\displaystyle\E\left(\int_\Sigma \Vert\hat{T}_n(x)-T(x)\Vert_d\right)$ go to zero, the same holds for $\hat{f}_n(x)$ provided that $f$ is $\alpha$-regular.
\item for almost all $x\in\Sigma$, for $a <\frac12$, $|\hat{g}_n(x)-g(x)| =O_\P(n^{-a})$ where $\hat{g}_n$ is either $\hat{f}_n$ or $\hat{T}_n$ and $g$ is either $f$ (provided it is $\alpha$-regular) or $T$.
\end{itemize}
\end{corointro}
In a first section, we state our hypothesis on the process and the dynamical systems. We also give the precise results.\\
The second section is devoted to the proofs.\\
In the last section we provide examples of dynamical systems satisfying our hypothesis. We also provide some  simulations. 
\tableofcontents
\section{Hypothesis and statement of the results}\label{def}
\subsection{Weak dependence}
As mentioned quickly in the introduction, we shall consider a class of weak mixing process with respect to a Banach space of bounded functions $\CC$. This functional definition of dependence is more general than strong mixing used in \cite{Mas} and \cite{FV}. We shall see that it  encounters a large class of dynamical systems (see also \cite{DP} for other examples than dynamical systems). This kind of functional definition has been introduced in \cite{DL} and \cite{DP} for Lipschitz or BV functions (see also \cite{taudep}). \\
For simplicity, let $\Sigma\subset \R^d$ and $(X_i)_{i\in\N}$ a process taking values in $\Sigma$. Our proofs could probably be extended to sets $\Sigma$ included in more general Banach spaces. In the following, $\Vert \ \Vert_d$ denotes the sup norm on $\R^d$~: $\Vert (x_1\/, \ldots\/, x_d)\Vert_d=\max(|x_i|\/, \ i=1\/, \ldots \/, d)$. 
\begin{defi}\label{weakmixing}
Let $\CC$ be a Banach space of bounded functions on $\Sigma$. We consider a norm on $\CC$ of the form~:
$$\Vert g \Vert_\CC = \Vert g \Vert_\infty + C(g)$$
where $C(\cdot)$ is a semi-norm on $\CC$ and $\Vert \cdot \Vert_\infty$ is the sup. norm on $\CC$. Let $\CC_1$ be the ``semi-ball'' of functions $g\in\CC$ such that $C(g)\leq 1$. Let $\MM_i$ be the $\sigma$-algebra generated by $X_0\/, \ldots\/, X_i$. The $\CC$-mixing coefficients are~:
\begin{align} 
\Phi_{\CC}(n) =& \sup \{ |\E(Z\/g(X_{i+n}))-\E(Z)\E(g(X_{i+n}))| \ i \in \N \/, \  Z \ \mbox{is} \nonumber\\
 & \ \MM_i\/-\/\mbox{measurable and} \ \Vert Z \Vert_1\leq 1\/,  \ g \in \CC_1 \}\/. \tag{*}\label{phi}
\end{align}
The process $(X_i)_{i\in\N}$ is $\CC$-weakly dependant if $(\Phi_{\CC}(n))_{n\in\N}$ is summable. \\
Let $Y_i$ be a $L^1$ stationary process, taking values in $\R$.  Let $\widetilde{\MM}_i$ be the $\sigma$-algebra generated by $X_0\/,Y_0\/, \ldots\/, X_i\/,Y_i$. We shall say that $(X_i)_{i\in\N}$ is $\CC$-weakly dependant with respect to $(Y_i)_{i\in\N}$ if the coefficients~:
\begin{align}
\Phi_{\CC\/,(Y_i)_{i\in\N}}(n) =& \sup \{ |\E(Z\/Y_{i+n}\/g(X_{i+n}))-\E(Z)\E(Y_{i+n}\/g(X_{i+n}))| \ i \in \N \/, \  Z \ \mbox{is} \nonumber\\
 & \ \widetilde{\MM}_i\/-\/\mbox{measurable and} \ \Vert Z \Vert_1\leq 1\/,  \ g \in \CC_1 \}\/, \tag{**}\label{phi2}
\end{align}
are summable.
\end{defi} 
\begin{rem}
Examples of spaces $\CC$ that we shall consider are~: 
\begin{itemize}
\item the space of function with bounded variations, in that case, $C(g)$ is the total variation of $g$~; 
\item the space of Lipschitz (resp. H\"older) functions, in that case, $C(g)$ is the Lipschitz (resp. H\"older) constant~;
\item the space of $C^1$ function could also be considered, in that case, $C(g)$ is the sup norm of $g'$.
\end{itemize}
\end{rem}
\subsection{Regularity conditions}
We now state the conditions on the kernel $K$, the regression function $r$ and the density function $f$.
\begin{cond}\label{Kernel}
The kernel $K$ is a nonnegative function in $\CC$ (so it is bounded),  with compact support $D$ and  with integral $1$~:
$$\int_\Sigma K \/ dm =1\/.$$
For $h>0$ and $x\in \Sigma$, let $K_{h\/,x}(t)=K\left(\frac{x-t}{h}\right)$. We assume that there exists $0 \leq \beta<1$ such that $C(K_{h\/,x}) \leq \frac{C(K)}{h^\beta}$.
\end{cond}
\begin{cond}\label{almosthol}
Consider a function $g$ on $\Sigma$ and $0<\alpha\leq 1$. Let $B_g(u\/,h)$ be the set of points $x$ such that 
$$\sup_{d(x\/,y)<h} |g(x)-g(y)|>u^{\alpha} \/.$$
If $g$ is a map from $\Sigma$ to $\Sigma$, $B_g(u\/,h)$ is the set of points $x$ such that 
$$\sup_{d(x\/,y)<h} \Vert g(x)-g(y)\Vert_d >u^{\alpha} \/.$$
For any decreasing to zero sequences $(u_n)_{n\in\N}$ and $(h_n)_{n\in\N}$, let $A_n = B_g(u_n\/,h_n)$ and $B_N = \displaystyle\bigcap_{n\geq N}\bigcup_{p\geq n} A_p$.\\
We shall say that $g$ is  $\alpha$-regular if for any sequence $(u_n)_{n\in\N}$ decreasing to $0$, for any sequence $(h_n)_{n\in\N}$ with $h_n=o(u_n)$, there exist constants $D(g)\/, D'(g)>0$ and $\gamma\/, \gamma' >0$ such that $m(B_N) \leq  D(g) \/ h_N^\gamma$ and $m(A_N)\leq D'(g)\/ h_N^{\gamma'}$.
\end{cond}
All our results are proved under the assumption that $r$ and/or $f$ are $\alpha$-regular for some $0<\alpha\leq 1$.
\begin{rems}
If $g$ is $\alpha$-H\"older then it is  $\alpha$-regular with $m(B_g(u\/,h))=0$ for all $u\geq H(g)h$ with $H(g)$ the H\"older constant of $g$. \\
If $g$ is Lipschitz on $\Sigma\setminus\{x\}$ and discontinuous at $x$ then $ m(B_g(u\/,h))\leq \mbox{Ct}h^d$ for all $u\geq L(g) h$ where $L(g)$ is the Lipschitz constant of $g$ on $\Sigma\setminus\{x\}$.\\
Functions of bounded variations satisfy this condition (see section \ref{ex-simul}).\\
An  $\alpha$-regular function  need not to be H\"older  everywhere but we control the points where this is not the case.
\end{rems}
\subsection{Main result}
Let us state our main result for $\CC$-weak dependant process. 
\begin{theo}\label{theomain}
Assume that $(X_i)_{i\in\N}$ is a stationary process absolutely continuous and is $\CC$-weakly dependant with respect to $(Y_i)_{i\in\N}$, $(Y_i)_{i\in\N}$ is a stationary  process taking values in $\R$. Let $r$ be the regression function ($r=\E(Y_i|X_i)$) and $f$, the density function, we assume that  $\inf f>0$. Then there exists $M>0$, $L>0$, $L'>0$, $R>0$, a set $A_{u\/,h}$  such that  for all $n\in \N^\star$,  for all $t\in\R^{+}$,  for $x\not\in A_{u\/,h}$, for all $u\geq h\cdot \mbox{Diam}(D)$, $m(A_{u\/,h})\leq R h^{\gamma'}$,\\
if $f$ is $\alpha$-regular,
\begin{equation}\label{expof2}
\P(|\widehat{f}_n(x)-f(x)|>t-u^\alpha)\leq 2e^{\frac1e} \exp[-t^2 M n h^{\beta+2}]\/,
\end{equation}
if $r$ is $\alpha$-regular and bounded,
\begin{equation}\label{expor2}
\P(|\widehat{r}_n(x) -r(x)|>t-u^\alpha)\leq e^{\frac1e} \left(2\exp[-t^2 L n h^{\beta+2}]+ \exp[-L'\/n\/ h^{\beta+2}]\right) \/,
\end{equation}
if $r$ is $\alpha$-regular and $Y_i\in L^\infty$,
\begin{equation}\label{expor2borne}
\P(|\widehat{r}_n(x) -r(x)|>t-u^\alpha)\leq 2e^{\frac1e} \exp[-t^2 L n h^{\beta+2}] \/.
\end{equation}
As a consequence,  provided $h=h_n$ goes to zero and $n h_n^{\beta+2}=O(n^\eps)$, $\eps>0$, we obtain the following convergences~:
\begin{itemize}
\item for $m$ almost all $x\in \Sigma$,  $\hat{r}_n(x)$ converges to $r(x)$ almost surely and in $L^p$ for any $1\leq p$, the same holds for $\hat{f}_n(x)$ provided $f$ is $\alpha$-regular.
\item $\displaystyle\E\left(\int_\Sigma |\hat{r}_n(x)-r(x)|\right)$ goes to zero, the same holds for $\hat{f}_n(x)$ provided $f$ is $\alpha$-regular.
\item for almost all $x\in\Sigma$, for $a <\frac12$, $|\hat{r}_n(x)-r(x)| =O_\P(n^{-a})$, the same holds for $\hat{f}_n(x)$ provided $f$ is $\alpha$-regular.
\end{itemize}
\end{theo}
\subsection{Dynamical systems}\label{dyn}
We turn now to our main motivation~: dynamical systems. Consider a dynamical system $(\Sigma\/, T\/, \mu)$. That is, $T$ maps $\Sigma$ into itself, $\mu$ is a $T$-invariant probability measure on $\Sigma$, absolutely continuous with respect to the Lebesgue measure $m$. We assume that the dynamical system satisfy the following mixing property~: for all $\varphi\in L^1(\mu)$, $\psi\in \CC$,
\begin{equation}\label{mixing}
\left|\int\limits_\Sigma \psi \cdot \varphi\circ T^n \/ d\mu - \int\limits_\Sigma \psi \/ d\mu \int\limits_\Sigma \varphi \/ d\mu \right| \leq \Phi(n) \Vert \varphi \Vert_1 \/ \Vert \psi \Vert_{\CC} \/, 
\end{equation}
with $\Phi(n)$ summable.
\begin{defi}\label{goodclass}
A dynamical system $(\Sigma\/, T\/, \mu)$ is of the good class if 
\begin{enumerate}
\item $\Sigma$ is a bounded subset of $\R^d$, 
\item the map $T~:~\Sigma~\rightarrow~\Sigma$ is  $\alpha$-regular,
\item the invariant density $f$ ($\mu=f\/m$) verify $\inf f>0$,
\item it satisfy the mixing property (\ref{mixing}).
\end{enumerate}
\end{defi}
We shall also assume that the Banach space $\CC$ is {\em regular} in the sense that for any $\varphi\in\CC$, there exists $R(\varphi) \in \R$ such that 
\begin{equation}\label{regspace}
\Vert \varphi + R(\varphi)\Vert_\infty \leq C(\varphi) \/.
\end{equation}
We consider the stationary process taking values in $\Sigma$~: $X_i=T^i$ with law $\mu$ and $Y_i=X_{i+1}+\eps_i$ where the $\eps_i$ are $L^1$ independent random vectors, with independent coordinates, identically distributed random variables with $0$ mean and independent on $(X_i)_{i\in\N}$. Then $T$ is the regression function $T(x)=\E(Y_i\/|\/ X_i=x)$. Let $\hat{r}_n^j$ is the estimator for   $Y_i^j=X_{i+1}^j+\eps_i^j$, defined by (\ref{defreg}),  $j=1\/, \ldots \/, d$, $X_i^j$ is the $j$th coordinate of $X_i$ and $\eps_i^j$ is the $j$th coordiante of $\eps_i$. We shall denote $\hat{T}_n$  the estimator of $T$ which has coordinates $\hat{r}_n^j$. \\
Theorem \ref{theomain} applied in the context of dynamical systems gives the following result.
\begin{coro}\label{corodyn}
Let $T$ be a dynamical system of the good class, assume that the  Banach space $\CC$ is regular. There exists $M>0$, $L>0$, $R>0$, $0\leq \beta \leq 1$, $\gamma'>0$  such that outside a set of measure less than $R h^{\gamma'}$,   for all $t\in\R^{+}$,   for all $u\geq \mbox{Ct} \/h$, 
\begin{eqnarray*}
\P(\Vert\hat{T}_n(x)-T(x)\Vert_d>t-u^\alpha) &\leq &2\/d\/e^{\frac1e} \exp[-t^2 L h^{\beta+2}  n] \\
\P(|\hat{f}_n(x)-f(x)|>t-u^\alpha)&\leq &2e^{\frac1e} \exp[-t^2 M h^{\beta+2} n] \ \mbox{if} \ f \ \mbox{is} \ \alpha\mbox{-regular}\/.
\end{eqnarray*}
As a consequence,  provided $h=h_n$ goes to zero  and $n h_n^{\beta+2} =O(n^\eps)$, $\eps>0$, we obtain the following convergences~:
\begin{itemize}
\item for $m$ almost all $\hat{T}_n(x)$ converges to $T(x)$ $x\in \Sigma$, $\hat{f}_n(x)$ converges to $f(x)$ almost surely and in $L^p$  for any $1\leq p$, if $f$ is $\alpha$-regular, $\hat{f}_n(x)$ converges to $f(x)$ almost surely and in $L^p$  for any $1\leq p$,
\item $\displaystyle\E\left(\int_\Sigma \Vert\hat{T}_n(x)-T(x)\Vert_d\right)$  go to zero, if $f$ is $\alpha$-regular,\\
$\displaystyle\E\left(\int_\Sigma |\hat{f}_n(x) -f(x)| \/dx \right)$  go to zero.
\item for almost all $x\in\Sigma$, for $a<\frac12$, $|\hat{g}_n(x)-g(x)| =O_\P(n^{-a})$ where $\hat{g}_n$ is either $\hat{f}_n$ or $\hat{T}_n$ and $g$ is either $f$ (provided $f$ is $\alpha$-regular) or $T$.
\end{itemize}
\end{coro}
\section{Proofs}
Let us now prove the results stated in the above section. The main ingredient is an exponential inequality that has been  obtained in \cite{DP} in the case $\CC = BV$, and in \cite{taudep} for more general spaces $\CC$.
\subsection{Main ingredient~: exponential inequality}
\begin{prop} \label{expo1}
Let $(X_i)_{i\in\N}$ be a $\CC$-weakly mixing process. Let the coefficients $\Phi_{\CC}(k)$ be defined by (\ref{phi}). For $\varphi \in \CC$, $p\geq 2$, define  
$$S_n(\varphi) = \sum_{i=1}^n \varphi(X_i) $$
and
$$b_{i\/,n} = \Phi_\CC(0) \/\left(\sum_{k=0}^{n-i} \Phi_{\CC}(k) \right) \Vert \varphi(X_i) - \E(\varphi(X_i)) \Vert_{\frac{p}2}   C(\varphi) \/.$$ 
For any $p\geq 2$, we have the inequality~: 
\begin{eqnarray} 
\lefteqn{\Vert S_n(\varphi) - \E(S_n(\varphi)) \Vert_p \leq \left(2p \Phi_\CC(0)\sum_{i=1}^n b_{i\/,n} \right)^{\frac12}}\nonumber\\
&& \leq C(\varphi) \left(2p\sum_{k=0}^{n-1} (n-k)\Phi_{\CC}(k) \right)^{\frac12}\/. \label{hoeff_p}
\end{eqnarray} 
As a consequence, we obtain 
\begin{eqnarray}
\lefteqn{\P\left(|S_n(\varphi) - \E(S_n(\varphi))|>t \right)}\nonumber\\
&& \leq e^{\frac1e} \exp \left(\frac{-t^2}{2e (C(\varphi))^2 \Phi_\CC(0)\/\sum_{k=0}^{n-1} (n-k)\Phi_{\CC}(k)} \right)\/. \label{concen} 
\end{eqnarray} 
\end{prop} 
Applying Proposition \ref{expo1} to the function $\phi=K$ will provide an exponential inequality for $\hat{f}_n$. The proof of Proposition \ref{expo1} may be found in \cite{taudep} (Proposition 1.1), see also \cite{DP} for a version of this proposition with $\CC=BV$. In order to have an exponential inequality for $\hat{g}_n$, we need the following result. 
\begin{prop} \label{expo2}
Let $(X_i)_{i\in\N}$ be a $\CC$-weakly mixing process with respect to $(Y_i)_{i\in\N}$. Let the coefficients $\Phi_{\CC\/,(Y_i)_{i\in\N}}(k)$ be defined by (\ref{phi2}). For $\varphi \in \CC$, $p\geq 2$, define  
$$\widetilde{S}_n(\varphi) = \sum_{i=1}^n Y_i\/ \varphi(X_i) $$
and
$$\widetilde{b}_{i\/,n} = \Phi_{\CC(0)\/, (Y_i)_{i\in\N}}\left(\sum_{k=0}^{n-i} \Phi_{\CC\/,(Y_i)_{i\in\N}}(k) \right) \Vert Y_i\varphi(X_i) - \E(Y_i\varphi(X_i)) \Vert_{\frac{p}2}   C(\varphi) \/.$$ 
For any $p\geq 2$, we have the inequality~: 
\begin{eqnarray} 
\lefteqn{\Vert \widetilde{S}_n(\varphi) - \E(\widetilde{S}_n(\varphi)) \Vert_p \leq \left(2p \sum_{i=1}^n \widetilde{b}_{i\/,n} \right)^{\frac12}}\nonumber\\
&& \leq C(\varphi) \left(2p\/ \Phi_{\CC\/, (Y_i)_{i\in\N}}(0)\/  \sum_{k=0}^{n-1} (n-k)\Phi_{\CC\/,(Y_i)_{i\in\N}}(k) \right)^{\frac12}\/. \label{hoeff_p2}
\end{eqnarray} 
As a consequence, we obtain 
\begin{eqnarray}
\lefteqn{\P\left(|S_n(\varphi) - \E(S_n(\varphi))|>t \right)}\nonumber\\
&& \leq e^{\frac1e} \exp \left(\frac{-t^2}{2e (C(\varphi))^2 \/  \Phi_{\CC(0)\/, (Y_i)_{i\in\N}}\/ \sum_{k=0}^{n-1} (n-k)\Phi_{\CC\/,(Y_i)_{i\in\N}}(k)} \right)\/. \label{concen2} 
\end{eqnarray} 
\end{prop}  
\begin{proof}
Proposition \ref{expo2} is proved as Proposition 1.1 in \cite{taudep} using the variable $Z_i=Y_i\/\varphi(X_i) - \E(Y_i\/\varphi(X_i))$ for which we control $\Vert \E(Z_k|\widetilde{\MM}_i) \Vert_\infty$ with the coefficients $\Phi_{\CC\/,(Y_i)_{i\in\N}}(k-i)$, for $k\geq i$.
\end{proof}
\subsection{Proof of Theorem \ref{theomain}}
As already mentioned, the proof of Theorem \ref{theomain} uses Propositions \ref{expo1} and \ref{expo2} applied to $\varphi = K$. This gives a control on the deviation between $\hat{f}_n(x)$ (resp. $\hat{g_n}(x)$) and $\E(\hat{f}_n(x))$ (resp. $ \E(\hat{g}_n(x))$). Then it remains to  see to that   $\E(\hat{f}_n(x))$  is close to $f(x)$. In order to obtain the result for the regression function, we deduce from the previous results a control on the deviation between $\hat{r}_n(x)$ and $\displaystyle \frac{\E(\hat{g}_n(x))}{\E(\hat{f}_n(x))}$, then we prove that $\displaystyle \frac{\E(\hat{g}_n(x))}{\E(\hat{f}_n(x))}$ is close to $r(x)$. To this aim, we begin with two lemmas.
\begin{lem}\label{espf}
We assume that the density function $f$ is  $\alpha$-regular. For all $u\geq h \/ \mbox{diam} (D)=:k$, if $x\not\in B_f(u\/, k)$ we have~:
$$|\E(\hat{f_n}(x)) -f(x)| \leq u^\alpha \/.$$ 
\end{lem}
\begin{proof}
We have~:
$$\E(\hat{f_n}(x)) = \int_D K(y) f(x-hy) \/ dm(y) \/.$$
If $x\not\in B_f(u\/, k)$ then for $y \in D$, 
$$|f(x)-f(x-hy)| \leq u^\alpha \/.$$
The result follows from the fact that the integral of $K$ is $1$.
\end{proof}
\begin{lem}\label{espr}
We assume that the regression function $r$ is $\alpha$-regular. For all $u\geq h \/ \mbox{diam} (D)=:k$, if $x\not\in B_r(u\/, k)$ we have~: 
$$\left|\frac{\E(\hat{g}_n(x))}{\E(\hat{f_n}(x))} -r(x)\right| \leq u^\alpha \/.$$ 
\end{lem}
\begin{proof}
We have~:
$$\frac{\E(\hat{g}_n(x))}{\E(\hat{f_n}(x))} = \frac{\displaystyle \int_D r(x-yh) K(x) f(x-yh) \/ dm(y)}{\displaystyle\int_D K(x) f(x-yh) \/ dm(y)}\/.$$
If $x\not\in B_r(u\/, k)$ then for $y \in D$, 
$$|r(x)-r(x-hy)| \leq u^\alpha $$
and the result follows.
\end{proof}
\begin{proof}[Proof of Theorem \ref{theomain}]
Proposition \ref{expo1} applied to $\varphi(t)=K\left(\frac{x-t}{h}\right)$ gives~: for all $t>0$, for all $n\in\N$, for all $x\in\Sigma$,
$$\P(|\hat{f}_n(x) - \E(\hat{f}_n(x))|>t) \leq e^{\frac1e} \exp\left(-\frac{t^2 \/n\/h^{\beta+2}}{2e\/R\/ \Phi_\CC(0)\/ C(K)} \right)\/,$$
where $R$ is the smallest positive number such that 
$$\sum_{k=0}^{n-1}(n-k) \Phi_\CC(k) \leq R \/ n \/.$$
This together with Lemma \ref{espf} gives (\ref{expof2}) for $x\not\in B_f(u\/,k)$ ($k=h\/diam(D)$). \\
In order to obtain the estimation for the regression function $r$, we apply Proposition \ref{expo2} to $\varphi(t)=K\left(\frac{x-t}{h}\right)$. We obtain~:
$$\P(|\hat{g}_n(x) - \E(\hat{g}_n(x)|>t) \leq e^{\frac1e} \exp\left(-\frac{t^2 \/n\/h^{\beta+2}}{2e\/R'\/\Phi_{\CC\/,(Y_i)_{i\in\N}}(0)\/  C(K)} \right)\/,$$
where $R'$ is the smallest positive number such that 
$$\sum_{k=0}^{n-1}(n-k) \Phi_{\CC\/,(Y_i)_{i\in\N}}(k) \leq R' \/ n \/.$$
Now, for any $t>0$, 
\begin{eqnarray*}
\P\left(\left|\hat{r}_n(x) - \frac{\E(\hat{g}_n(x))}{\E(\hat{f_n}(x))}\right|>t \right) &\leq& \P\left(|\hat{g}_n(x) - \E(\hat{g}_n(x)|>\frac{t}2 \E(\hat{f}_n(x) \right) \\
&+& \P\left(|\hat{f}_n(x) - \E(\hat{f}_n(x)|>\frac{t}2 \E(\hat{f}_n(x) \hat{r}_n(x)^{-1}\right) \/.
\end{eqnarray*}
Now, assume that $Y_i\in L^{\infty}$, we remark that $|\hat{r}_n(x)| \leq \Vert Y_i\Vert_\infty =: y_{max}$ thus~:
\begin{eqnarray*}
\P\left(\left|\hat{r}_n(x) - \frac{\E(\hat{g}_n(x))}{\E(\hat{f_n}(x))}\right|>t \right)&\leq& \P\left(|\hat{g}_n(x) - \E(\hat{g}_n(x))|>\frac{t}2 \inf f \right)\\
&+&  \P\left(|\hat{f}_n(x) - \E(\hat{f}_n(x))|>\frac{t}2 \frac{\inf f}{y_{max}}\right) \/.
\end{eqnarray*}
Propositions \ref{expo1} and \ref{expo2} give~:
$$\P\left(\left|\hat{r}_n(x) - \frac{\E(\hat{g}_n(x))}{\E(\hat{f_n}(x))}\right|>t \right)\leq 2 e^{\frac1e} \exp[-C t^2 h^{\beta+2}\/n]$$
where 
$$C = \max\left(\frac{(\inf f)^2}{8e\/R'\/\Phi_{\CC\/,(Y_i)_{i\in\N}}(0)\/ C(K)} \/, \frac{(\inf f)^2}{8e\/R\/\Phi_{\CC}(0)\/ C(K)\/ y^2_{\max}} \right) \/.$$
Lemma \ref{espr} gives  (\ref{expor2borne}) for $x\not\in B_r(u\/, k)$. Let $A_h= B_f(u\/, k)\cup  B_r(u\/, k)$, the first part of Theorem \ref{theomain} is now proved for $x\not\in A_h$ if $Y_i\in L^{\infty}$.\\
If we don't assume $Y_i\in L^{\infty}$, but $r$ is bounded by $r_{\max}$, we write~:
\begin{eqnarray*}
\lefteqn{\P\left(\left|\hat{r}_n(x) - \frac{\E(\hat{g}_n(x))}{\E(\hat{f_n}(x))}\right|>t \right)\leq  \P\left(|\hat{g}_n(x) - \E(\hat{g}_n(x))|>\frac{t}2 \hat{f}_n(x) \right)}\\
&&+\P\left(|\hat{f}_n(x) - \E(\hat{f}_n(x))|>\frac{t}2 \hat{f}_n(x) \left[\frac{\E(\hat{g}_n(x))}{\E(\hat{f_n}(x))}\right]^{-1}\right)\\
&\leq & \P\left(|\hat{g}_n(x) - \E(\hat{g}_n(x))|>\frac{t}4 \inf f \right)+\P\left(|\hat{f}_n(x) - \E(\hat{f}_n(x))|>\frac{t}4 \frac{\inf f}{r_{\max}}\right) \\
& & + \P \left( \hat{f_n}(x)< \frac{\inf f}2\right) \/.
\end{eqnarray*}
Propositions \ref{expo1} and \ref{expo2} give~:
$$\P\left(|\hat{r}_n(x) - \frac{\E(\hat{g}_n(x))}{\E(\hat{f_n}(x))}|>t \right)\leq e^{\frac1e} \left[2\/\exp[-C't^2h^{\beta+2}\/n] + \exp[-C''\/n\/h^{\beta+2}] \right]\/,$$
where $C'$ and $C''$ may be expressed in terms of $\inf f$, $r_{\max}$, $R$, $R'$,$\Phi_{\CC}(0)$, $\Phi_{\CC\/,(Y_i)_{i\in\N}}(0)$, $C(K)$. Then (\ref{expor2}) follows as above for $x\not\in A_h$.\\
Let us conclude with the almost everywhere and $L^p$ convergence, for all $1\leq p$. Let us begin with the convergence in $L^p(m)$. We fix $1\leq p$. We remark that since the kernel $K$ is bounded, so is $\hat{f}_n$~: $\sup(\hat{f}_n) \leq \frac{\sup K}{h}$. Also, if $x\not\in B_f(u\/,\eps)$ then $f(x)\leq \frac{1+u^\alpha}{\eps^d}$. Let $E_{t\/,u}(x)$ be the event
$$E_{t\/,u}(x)= \{|\hat{f}_n(x) -f(x)|>t-u^\alpha \} \/,$$
we have~:
\begin{eqnarray*}
\int |\hat{f}_n(x) -f(x)|^p \/ d\P & = & \int_{E_{t\/,u}(x)}  |\hat{f}_n(x) -f(x)|^p \/ d\P \\
&&+ \int_{[E_{t\/,u}(x)]^c}  |\hat{f}_n(x) -f(x)|^p \/ d\P\\
&\leq & \P(E_{t\/,u}(x)) \left(\frac{1+u^\alpha}{diam(D) h}+ \frac{\sup K}{h} \right)^p + (t-u^\alpha)^p \\
&\leq& \frac{\mbox{Ct}}{h^p}  \exp[-C t^2 h^\beta\/n] + (t-u^\alpha)^p\/.
\end{eqnarray*}
Take $u$ such that $u^\alpha=\frac1{\ln n}$, $t=2\/ u^\alpha$ and $h = h_n= O(\frac1{n^\xi})$ with $0<\xi< \frac1{\beta+2}$. Then for $x\not\in A_n=B_f(u_n\/,h_n)$, there exists $0<\kappa<1$ and constants $C_1\/,C_2>0$ such that~:
$$\Vert \hat{f}_n(x) -f(x)\Vert_p^p \leq C_1 e^{-C_2n^\kappa}+ \left(\frac1{\ln n}\right)^p\/.$$
Thus, $\Vert \hat{f}_n(x) -f(x)\Vert_p$ goes to zero for $x\not\in B_N=\cap_{n\geq N}\cup_{p\geq n} A_p$ and $m(B_N)=O(h_N^\gamma)$. We conclude that for almost all $x\in\Sigma$, $\hat{f}_n(x)$ goes to $f(x)$ in $L^p(m)$.  \\
The proof of the almost everywhere convergence is in the same spirit. Consider a sequence $t_m$ decreasing to $0$, for $x\in \Sigma$, 
$$\P(\{\hat{f}_n(x) \not\rightarrow f(x) \}) \leq \lim_{m\rightarrow \infty} \lim_{N\rightarrow \infty} \sum_{n\geq N} \P(\{|\hat{f}_n(x)-f(x)|>t_m\}) \/.$$
Take $h=h_n= O(\frac1{n^\xi})$ with $0<\xi< \frac1{\beta+2}$ then 
$$\lim_{N\rightarrow \infty} \sum_{n\geq N} \P(\{|\hat{f}_n(x)-f(x)|>t_m\}) = 0 $$
provided $t_m> u^\alpha$ with $u^\alpha > \mbox{Ct} h_n$ and $x\not\in B_{u\/,h_n}$, we take $u=u_n=\frac1{\ln n}$ and $A_n= B_{u_n\/, h_n}$, for $x\not\in \displaystyle\cap_{N}\cup_{n\geq N}   A_n =: B$, $\hat{f}_n(x)$ goes to $f(x)$  and $m(B)=0$.\\
The proofs of $L^p(m)$ and almost everywhere convergence for $\hat{r}_n$ follow the same lines.\\
Finally, to prove the bound for the convergence in probability, recall that $u_n=O_{\P}(n^{-a})$ if and only if
$$\lim_{M\rightarrow\infty}\limsup_{n\rightarrow\infty}  \P(n^a |u_n| >M) =0 \/.$$
Following the same lines as above, we prove that for $g$ either $f$ or $r$ and $\hat{g}_n$ either $\hat{f}_n$ or $\hat{r}_n$,
$$\lim_{M\rightarrow\infty}\limsup_{n\rightarrow\infty}  \P(n^a |\hat{g}_n(x)-g(x)| >M) =0$$
for $a<\frac12$, $x\not\in B=\cap_{n\geq 0}\cup_{p\geq n} A_p$ and $m(B)=0$.
\end{proof}
Let us show how we can apply the above result to dynamical systems.
\subsection{Dynamical systems and time reversed process}
It turns out that, in general, the process $(X_i)_{i\in\N}$ associated to a dynamical system $(\Sigma\/, T\/, \mu)$ is not weakly dependent. Nevertheless the condition of weak dependence is satisfied for a ``time reversed process'' whose law is the same as $(X_i)_{i\in\N}$. Indeed, Condition \ref{mixing} together with the fact that the space $\CC$ is regular (recall (\ref{regspace})) gives  (see \cite{DP} or \cite{taudep} for the details)~: for all $i\in \N$, for  $\psi\in \CC_1$, $\varphi\in L^1$, $\Vert \varphi\Vert_1\leq 1$,
\begin{equation}
| \mbox{Cov}(\psi(X_i)\/, \varphi(X_{n+i})) |\leq 2\Phi(n) \/.\label{dynmix}
\end{equation}
So that, if we consider the process $(\widetilde{X_i})_{i\in\N}$ defined by 
$$(\widetilde{X_0}, \ldots \/, \widetilde{X_n}) \stackrel{\mbox{Law}}{=} (X_n\/, \ldots \/, X_0) \ \forall n\in\N \/,$$
it is $\CC$-weakly dependent.\\
Recall that $Y_i=X_{i+1}+\eps_i$ with $(\eps_i)_{i\in\N}$ independent random vectors, with independent coordinates,  of zero mean, which is also independent of $(X_i)_{i\in\N}$.\\
In order to estimate the regression function $\E(Y_i | X_i)=T$, we shall estimate separately each coordinates of $X_{i+1}$. Let
\begin{eqnarray*}
\widehat{g}_n^j(x) &=& \frac1{nh} \sum_{i=0}^{n-1} Y_i^j K\left(\frac{x-X_i}h\right) \\
& =& \frac1{nh}\sum_{i=0}^{n-1} X_{i+1}^j K\left(\frac{x-X_i}h\right) +   \frac1{nh} \sum_{i=0}^{n-1} \eps_i^j  K\left(\frac{x-X_i}h\right)\/.
\end{eqnarray*}
Equation (\ref{dynmix}) implies also that  $(\widetilde{X}_i)_{i\in\N}$ is  $\CC$-weakly dependant with respect to each  $(\widetilde{X}_{i+1}^j)_{i\in\N}$, and also to each $(\eps_i^j)$, $j=1\/,\ldots \/, d$. We obtain inequalities (\ref{expor2borne}) separately for the two parts of $\widehat{g}_n^j(x)$ and then the two parts of $\widehat{r}_n^j(x)$. Then we deduce Corollary \ref{corodyn}.
\section{Examples and simulations}\label{ex-simul}
We shall now give some examples of discrete dynamical systems satisfying (\ref{dynmix}), such that $T$ is  $\alpha$-regular, admits a unique invariant density (with respect to the Lebesgue measure), $f$ which may be  also  $\alpha$-regular. \\
There is a large class of dynamical systems satisfying (\ref{dynmix}), we refer to the literature on dynamical systems~: \cite{Y1,Y2,tower,Li,Co,BuM1,BuM2,LSV,Br} and many other, see \cite{Ba} for a review on these topics. Below we consider Lasota-Yorke maps, unimodal maps, piecewise expanding maps in higher dimension.  Our results should also apply for hyperbolic maps but a more intricate study on the invariant density is required. 
\subsection{In dimension one~: piecewise expanding maps}
We shall consider piecewise expanding maps or ``Lasota-Yorke'' maps. Consider  $I=[0\/,1]$,  partitioned into subintervals $I_j$. On the interior of the $I_j$'s, the map $T$ is $C^2$, uniformly expanding $|T'|\geq \lambda>1$, and continuous on the closure of the $I_j$'s. Under additional conditions of mixing or covering (see below), it is well known (\cite{Br,LSV,Li,Co}) that $T$ admits a unique absolutely continuous invariant measure whose density $f$ belongs to the space $BV$ of functions of bounded variations. \\
Before we give precisely the hypothesis on the map, we prove that functions of bounded variations are  regular. We shall denote $m$ the Lebesgue measure on $I$. \\
Recall that a function $g$ on $I$ belongs to $BV$ if
$$\bigvee g = \sup \sum_{i=0}^n |g(a_i)-g(a_{i+1}) | < \infty $$
where the supremum is taken over all the finite subdivisions of $I$~; $a_0=0 < a_1 <\cdots < a_{n+1}=1$, then $\bigvee g$ is the total variation of $g$.
\begin{lem}\label{BV}
Let $g$ belongs to $BV$. Then for any sequences $(u_n)_{n\in\N}$ and $(h_n)_{n\in\N}$ decreasing to zero with $h_n^{2-\alpha}\leq u_n^2$, $0<\alpha<2$, let $B_g(u\/,h)$ be the set of points $x$ such that 
$$\sup_{d(x\/,y)<h} |g(x)-g(y)|>u\/,$$
$$A_n=  B_g(u_n\/,h_n) \ \mbox{and} \ \   B_N = \displaystyle\bigcap_{n\geq N}\bigcup_{p\geq n} A_p \/.$$
Then
$$m(B_N) \leq 3(\bigvee g) h_N^{\frac{\alpha}2}  \ \mbox{and} \ \ m(A_n)\leq 2 (\bigvee g) \frac{h_n}{u_n} $$
\end{lem}
\begin{proof}
Let $x_0=\inf B_N$ and $x_1$ such that $0\leq x_1-x_0\leq h_N$ and $x_1\in B_N$. Then $x_1$ belongs to $A_p$ for infinitely many $p\geq N$. In particular, there exists $p_1\geq N$ and $x_2$ with $|x_1-x_2|\leq h_{p_1}$ and $|g(x_1)-g(x_2)|\geq u_{p_1}$. Set $x_0=x_{1\/,0}$, $x_1=x_{1\/,1}$ and $x_2=x_{1\/,2}$. We construct   sequences (maybe finite) $(x_{i\/,0}\/,x_{i\/,1}\/,x_{i\/,2})$  and $(p_i)$ such that~:
\begin{enumerate}
\item $p_i\geq N$ for all $i\geq 1$,
\item $x_{i\/,0}\leq \min(x_{i\/,1}\/, x_{i\/,2})$~; $x_{i+1\/,0}\geq \max(x_{i\/,1}\/, x_{i\/,2})$, $x_{i\/,1} \in A_{p_i}$,
\item $|x_{i\/,0}-x_{i\/,1}|\leq h_{p_{i-1}}$ for all $i\geq 1$, with the convention that $h_{p_0}=h_N$,
\item $|x_{i\/,1}-x_{i\/,2}|\leq h_{p_i}$ for all $i\geq 1$,
\item $|g(x_{i\/,1})-g(x_{i\/,2})| \geq u_{p_i}$ for all $i\geq 1$,
\item $B_N \subset \bigcup_{i\geq 1} J_i$ with $J_i = [a_i\/,b_i+h_{p_{i-1}}]$, $ [a_i\/,b_i]= [x_{i\/,0}\/,x_{i\/,1}]\cup [x_{i\/,1}\/,x_{i\/,2}]$.
\end{enumerate}
If $x_{i\/,0}\/,x_{i\/,1}\/,x_{i\/,2}$ are already constructed, let $x_{i+1\/,0} = \inf [x\in B_N \/, \ x \geq \max(x_{i\/,1}\/, x_{i\/,2})+ h_{p_i}]$. Now, there exists $x_{i+1\/,1} \in B_N$ with $0\leq x_{i+1\/,1}-x_{i+1\/,0}\leq h_{p_i}$. Since  $x_{i+1\/,1} \in B_N$, there exists $p_{i+1}\geq p_i$ and $x_{i+1\/,2}$ with $|x_{i+1\/,1}-x_{i+1\/,2}|\leq h_{p_{i+1}}$ and  $|g(x_{i+1\/,1})-g(x_{i+1\/,2})| \geq u_{p_{i+1}}$. In that way, all the points of $B_N$ smaller than $\max(x_{i+1\/,2}\/,x_{i+1\/,1})+h_{p_{i+1}}$ are in $\displaystyle \bigcup_{j=1}^{i+1}J_j$. The construction stops when $\{x\in B_N \/, \ x \geq \max(x_{i\/,1}\/, x_{i\/,2})+ h_{p_i}\}$ is empty.\\
As a consequence of this construction, we get~:
$$\sum_{i\geq 1} u_{p_i} \leq \bigvee g \ \mbox{remark that the} \ J_i\mbox{'s are disjoint} \/,$$
$$m(B_N)\leq 3\sum_{i\geq0} h_{p_i} \/.$$
Now, using Cauchy-Schwartz inequality,
\begin{eqnarray*}
\sum_{i\geq0} h_{p_i}& \leq& \left(\sum_{i\geq 1} \frac{h_{p_i}^2}{u_{p_i}} \right)^{\frac12}\cdot  \left(\sum_{i\geq1} u_{p_i} \right)^{\frac12}\\
&\leq & \left(\bigvee g\right)^{\frac12} \cdot h_N^{\frac\alpha2} \left(\sum_{i\geq1} \frac{h_{p_i}^{2-\alpha}}{u_{p_i}} \right)^{\frac12} \\
&\leq & \bigvee g \cdot h_N^{\frac\alpha2} \ \mbox{recall that} \ h_n^{2-\alpha}\leq u_n^2 \/.
\end{eqnarray*}
Finally, we have proven~: $m(B_N) \leq 3(\bigvee g) h_N^{\frac{\alpha}2}$.\\
The estimation on the measure of $A_n$ is much simpler. We may cover $A_n$ by balls of radius $2h_n$, centered on points $x$ such that $\sup_{d(x\/,y)<h_n} |g(x)-g(y)|>u_n$. Then, $m(A_n)\leq 2\cdot k\cdot h_n$ where $k$ is the number of such balls. Now, $\bigvee g \geq k \cdot u_n$, thus $k\leq \displaystyle \frac{\bigvee g}{u_n}$. 
\end{proof}
\subsubsection{Finite number of pieces}
This subsection is devoted to Lasota-Yorke maps with a finite number of intervals of monotonicity. The main differences between infinite and finite pieces are some additional technical assumptions in the infinite case.\\
We consider the following conditions that may be found in \cite{Co,Li}.
\begin{cond} \label{expanding}\ \\
\begin{enumerate} 
\item  The partition $(I_i)_{i=1\/,\ldots\/r}$ of $I$ is finite, denote by $\PP_k$ the partition of invertibility of $T^k$.
\item $T$  satisfies the covering property ~: for all $k\in \N$, there exists $N(k)$ such that for all $P\in \PP_k$,
$$T^{N(k)}P = [0\/,1]\/.$$   
\item On each $\overline{I_i}$, the map $T$ is $C^2$.
\end{enumerate}
\end{cond}
The following result is standard and may be found in (\cite{Co,Li,Ba}).
\begin{theo}
Under Assumptions \ref{expanding}, the map $T$ admits a unique absolutely invariant probability measure, its density $f$ belongs to $BV$ and $\inf f >0$. Moreover, if $\mu = f m$ is this invariant probability measure, then (\ref{mixing}) is satisfied with $\CC = BV$ and $\Phi(n)= \gamma^n$, $0<\gamma<1$~: there exists $0<\gamma<1$, $C>0$, such that for any $\psi\in BV$ and $\varphi \in L^1(\mu)$, for any $n\in\N$, 
$$\left|\int\limits_I \psi \cdot \varphi\circ T^n \/ d\mu - \int\limits_I \psi \/ d\mu \int\limits_I \varphi \/ d\mu \right| \leq C \/ \gamma^n\/ \Vert \varphi \Vert_1 \/ \Vert \psi \Vert_{BV}\/.$$
\end{theo}
We already know that $f$ is  regular by Lemma \ref{BV}. In order to have the conditions of Corollary \ref{corodyn}, it remains to prove that $T$ is regular. This is indeed clear because $T$ is piecewise $C^2$, with finitely many points of discontinuity. As a consequence, for any sequences $(u_n)_{n\in\N}$ and $(h_n)_{n\in\N}$ decreasing to zero with $h_N\leq \frac{u_N}{\sup|T'|}$,
$$m(A_N) \leq \sum_{i=1}^p |g_i| \/ \cdot h_N$$
where we have denoted by  $x_i$, $i=1\/,\ldots \/, p$ the points of discontinuity, and $g_i$ the gap at $x_i$~: $g_i = |T(x_i^-)-T(x_i^+)|$. Moreover, as soon as $u_N\leq \inf g_i$, we have that $A_n\supset A_{n+1}$, so $B_N=A_N$ and we have the required control on the measure of $A_N$. \\
Thus, Corollary \ref{corodyn} applies and we get the following result.
\begin{coro}\label{finitepieces}
Let $T$ satisfies Assumption \ref{expanding}, let $K$ be a Kernel belonging to $BV$, let $\hat{f}_N$ and $\hat{T}_N$ be the estimators of $f$ and $T$.  For all $0<\alpha<1$, here exists $M>0$, $L>0$, $R>0$,   such that outside a set of measure less than $R h^{\alpha}$,   for all $t\in\R^{+}$,   for all $u\geq \mbox{Ct} \/h$, 
\begin{align}
\P(|\hat{f}_n(x)-f(x)|>t-u^\alpha)\leq 2e^{\frac1e} \exp[-t^2 M n h^2] \\
\P(|\hat{T}_n(x)-T(x)|>t-u^\alpha) \leq 2e^{\frac1e} \exp[-t^2 L  n h^2] \/.
\end{align}
As a consequence,  provided $h=h_n$ goes to zero  and $n h_n^2 =O(n^\eps)$, $\eps>0$, we obtain the following convergences~:
\begin{itemize}
\item for $m$ almost all $x\in \Sigma$, $\hat{f}_n(x)$ converges to $f(x)$ and $\hat{T}_n(x)$ converges to $T(x)$ almost surely and in $L^p$  for any $1\leq p$,
\item $\displaystyle\E\left(\int_I |\hat{f}_n(x) -f(x)| \/dx \right)$ and $\displaystyle\E\left(\int_I |\hat{T}_n(x)-T(x)|\right)$ go to zero.
\item for almost all $x\in\Sigma$, for $a <\frac12$, $|\hat{g}_n(x)-g(x)| =O_\P(n^{-a})$ where $\hat{g}_n$ is either $\hat{f}_n$ or $\hat{T}_n$ and $g$ is either $f$ or $T$.
\end{itemize}
\end{coro}
\begin{proof}
We apply Corollary \ref{corodyn}. Since the Kernel $K$ belongs to $BV$, we may take $\beta=0$. Also, with our hypothesis, $T$ and $f$ are in $BV$ so they are regular according to Lemma \ref{BV}. 
\end{proof}
{\underline{Simulations for $\beta$-maps}}. For $\beta >1$, the $\beta$-map is $T_\beta(x)=\beta x \ \mbox{mod}\ 1$. This map has been widely studied. It is well known that it satisfies Assumptions \ref{finitepieces}. It is the simplest example of a non-markov map on the interval. There is an explicit formula for the invariant density (see \cite{T})~: 
$$f_\beta(x)= C \sum_{i\geq 0} \beta^{-(i+1)} \Id_{[0\/,T^i1](x)}\/,$$
where $C$ is a constant chosen so that $f_\beta$ has integral $1$.\\
We have performed some simulations for several values of of $\beta$, different Kernels and for several noises $\eps_i$.  These simulations are summarized in the following table. In this table, AMEf means absolute mean error for the density i.e.
$$\mbox{AMEf} = \frac1p \sum_{k=1}^p |\widehat{f}_n(x_k) - f(x_k)| \/,$$
where the $x_k$ are $p$ points regularly espaced on $I$ on which we have calculate $\widehat{f}_n$. Below, $p=200$. AMET means  absolute mean error for the map $T$ (same formulae). For the Kernel, $P_2$ is the degree $2$ polynomial~: $\frac34(1-x^2)$.
\begin{center}
\begin{tabular}{|c|c|c|c|c|c|c|} \hline
\multicolumn{1}{|c|}{n}&\multicolumn{1}{c|}{h}&\multicolumn{1}{c|}{AMEf}&\multicolumn{1}{c|}{AMET}&\multicolumn{1}{c|}{$\beta$}&\multicolumn{1}{c|}{kernel}&\multicolumn{1}{c|}{noise, $\eps_i$}\\ \hline
$10^4$~~~&$0.01$~~~&$0.08234419$~~&$0.008309136$~&$\frac{27}{11}$~~~&$P_2$&no, $\eps_i=0$\\ \hline
$10^4$~~~~&$0.005$~~~&$0.09906515$~~&$0.004301326$~&$\frac{27}{11}$~~~&$P_2$&no, $\eps_i=0$\\ \hline
$5\cdot 10^4$~~~&$0.007$~~&$0.04428149$~~&$0.005530895$~&$\frac{27}{11}$~~&$P_2$&no, $\eps_i=0$\\ \hline
$2\cdot 10^5$~~&$0.001$~~&$0.05107575$~~&$0.001799785$~&$\frac{27}{11}$~~~&$P_2$&no, $\eps_i=0$\\ \hline
$5\cdot 10^4$~~~&$0.007$~~&$0.05492035$~~&$0.003809815$~&$\frac{27}{11}$~~~&$\Id_{[-\frac12\/, \frac12]}$&no, $\eps_i=0$\\ \hline
$5\cdot 10^4$~~~&$0.007$~~~&$0.04728425$~~&$0.008303824$~&$\frac{27}{11}$~~&$P_2$&${\mathcal U}[-0.3\/,0.3]$\\ \hline
$2\cdot 10^5$~~~&$0.0005$~~~&$0.07806642$~~&$0.011328519$~&$\frac{27}{11}$~~&$P_2$&${\mathcal U}[-0.3\/,0.3]$\\ \hline
$10^4$~~&$0.01$~~&$0.07473928$~~&$0.020744986$~&$\frac{27}{11}$~~~&$P_2$&${\mathcal N}(0\/,0.3)$~~~~\\ \hline
$5\cdot 10^4$~~~&$0.007$~~&$0.04269281$~~&$0.011423138$~&$\frac{27}{11}$~~~&$P_2$&${\mathcal N}(0\/,0.3)$~~~~\\ \hline
$2\cdot 10^5$~~~&$0.001$~~~&$0.05107575$~~&$0.001799785$~&$\frac{27}{11}$~~~&$P_2$&${\mathcal N}(0\/,0.3)$~~~~\\ \hline
$5\cdot 10^4$~~~~&$0.007$~~~&$0.05329131$~~&$0.007722570$~&$\frac{27}{11}$~~~&$\Id_{[-\frac12\/, \frac12]}$& ${\mathcal U}[-0.3\/,0.3]$\\ \hline
\hline \hline
$10^4$~~~&$0.01$~~~&$0.08165332$~~&$1.648713e-02$&$\frac{46}{11}$~~~&$P_2$&no, $\eps_i=0$\\ \hline
$5\cdot 10^4$~~~~&$0.007$~~~~&$0.04259507$~~&$1.071092e-02$&$\frac{46}{11}$~~~&$P_2$&no, $\eps_i=0$\\ \hline
$2\cdot 10^5$~~~~&$0.001$~~~&$0.05249840$~~&$1.536396e-04$&$\frac{46}{11}$~~~~&$P_2$&no, $\eps_i=0$\\ \hline
$5\cdot 10^4$~~~&$0.007$~~~&$0.03810643$~~&$1.482175e-02$&$\frac{46}{11}$~~~&$P_2$&${\mathcal U}[-0.3\/,0.3]$\\\hline  
$5\cdot 10^4$~~~&$0.007$~~~&$0.03961733$~~&$1.763502e-02$&$\frac{46}{11}$~~~&$P_2$&${\mathcal N}(0\/,0.3)$~~~~\\\hline  
$5\cdot 10^4$~~~&$0.007$~~~&$0.05467709$~~&$7.079913e-03$&$\frac{46}{11}$~~~&$\Id_{[-\frac12\/, \frac12]}$&no\\  \hline
$5\cdot 10^4$~~~&$0.007$~~~&$0.05109682$~~&$1.036748e-02$&$\frac{46}{11}$~~~&$\Id_{[-\frac12\/, \frac12]}$&${\mathcal U}[-0.3\/,0.3]$\\  
\hline
\end{tabular}
\end{center}
\begin{rem}
It seems that there no real influence of the kind of Kernel nor on the noise. What seems more surprising is the fact that the estimators for $T$ are much better than for $f$. Also, remark that it seems that the best $h_n$ is not the same for $f$ and for $T$.
\end{rem}
The graphics below correspond to $\beta=\frac{27}{11}$, $K=\Id_{[-\frac12\/, \frac12]}$, $\eps_i \leadsto {\mathcal U}[-0.2\/,0.2]$, $n=50000$ and $h=0.007$.\\
\begin{center}
\includegraphics[scale=0.43]{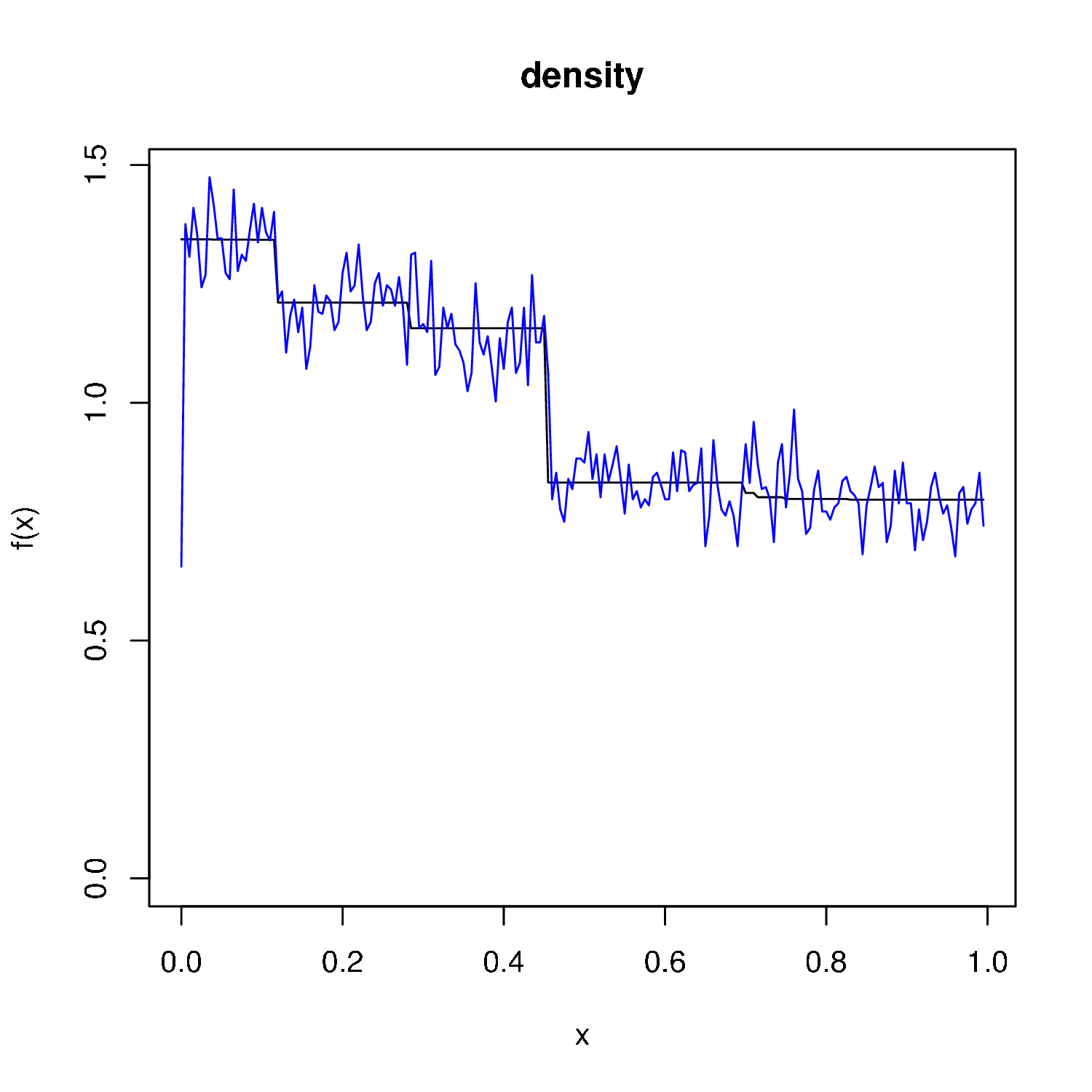}
\includegraphics[scale=0.43]{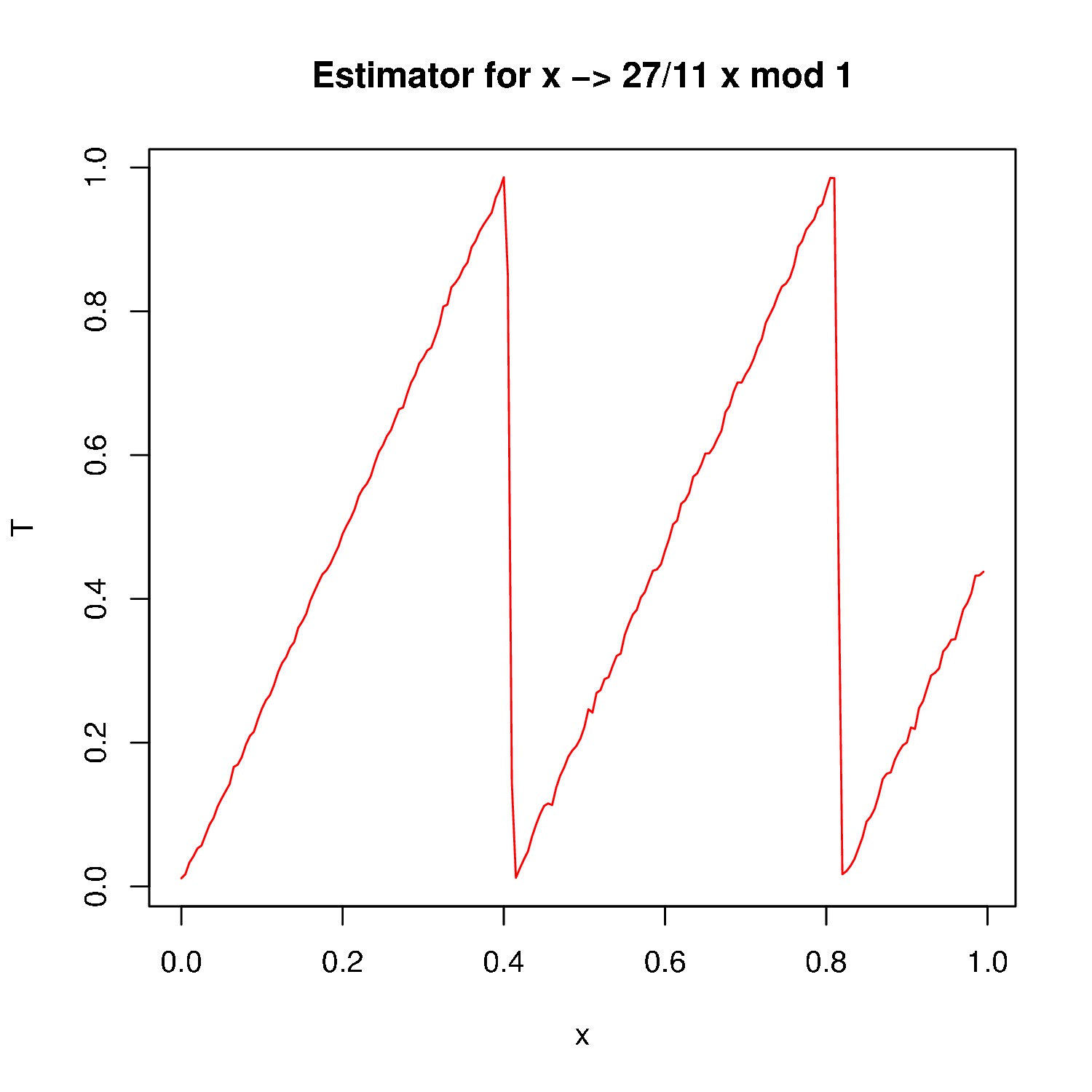}
\includegraphics[scale=0.43]{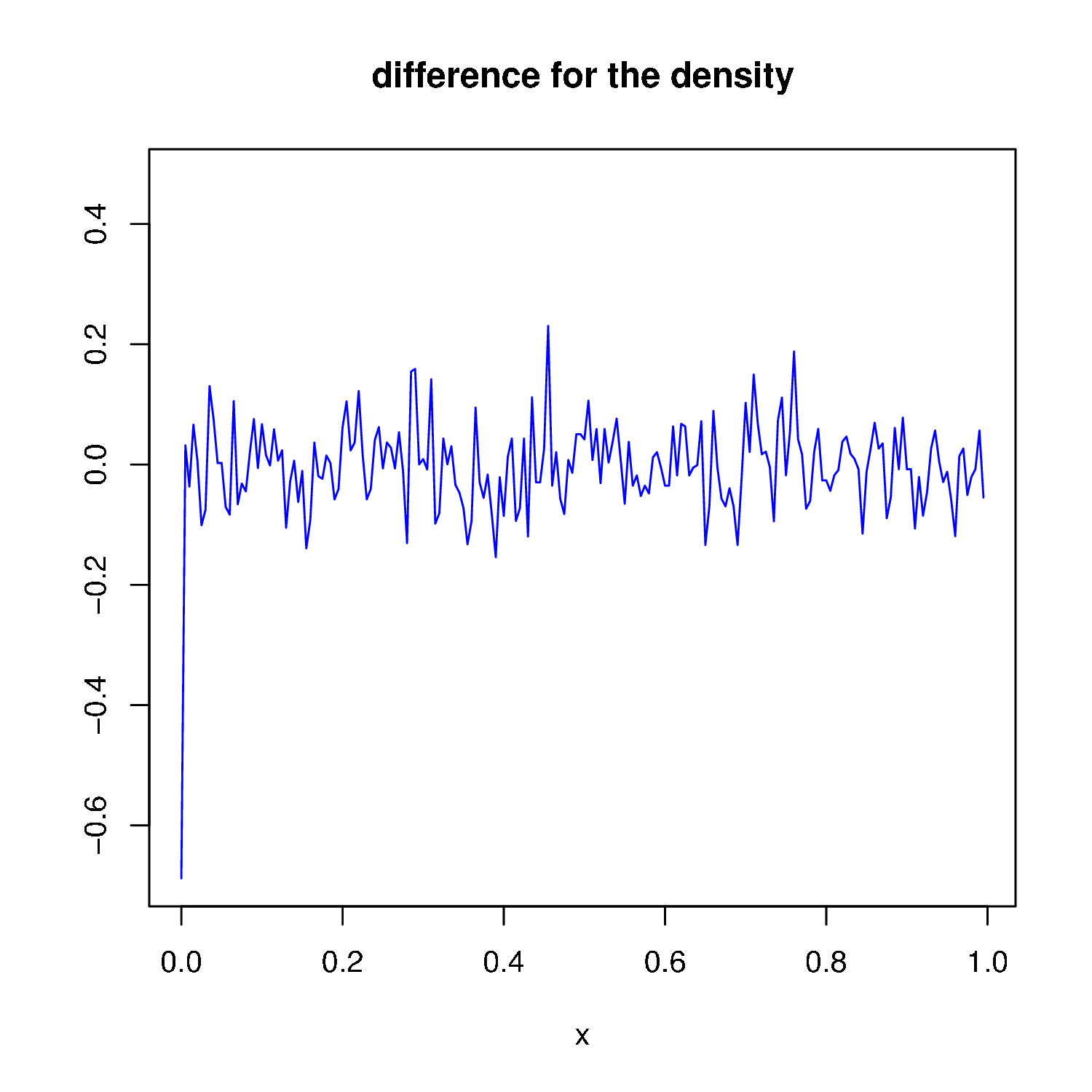}
\includegraphics[scale=0.43]{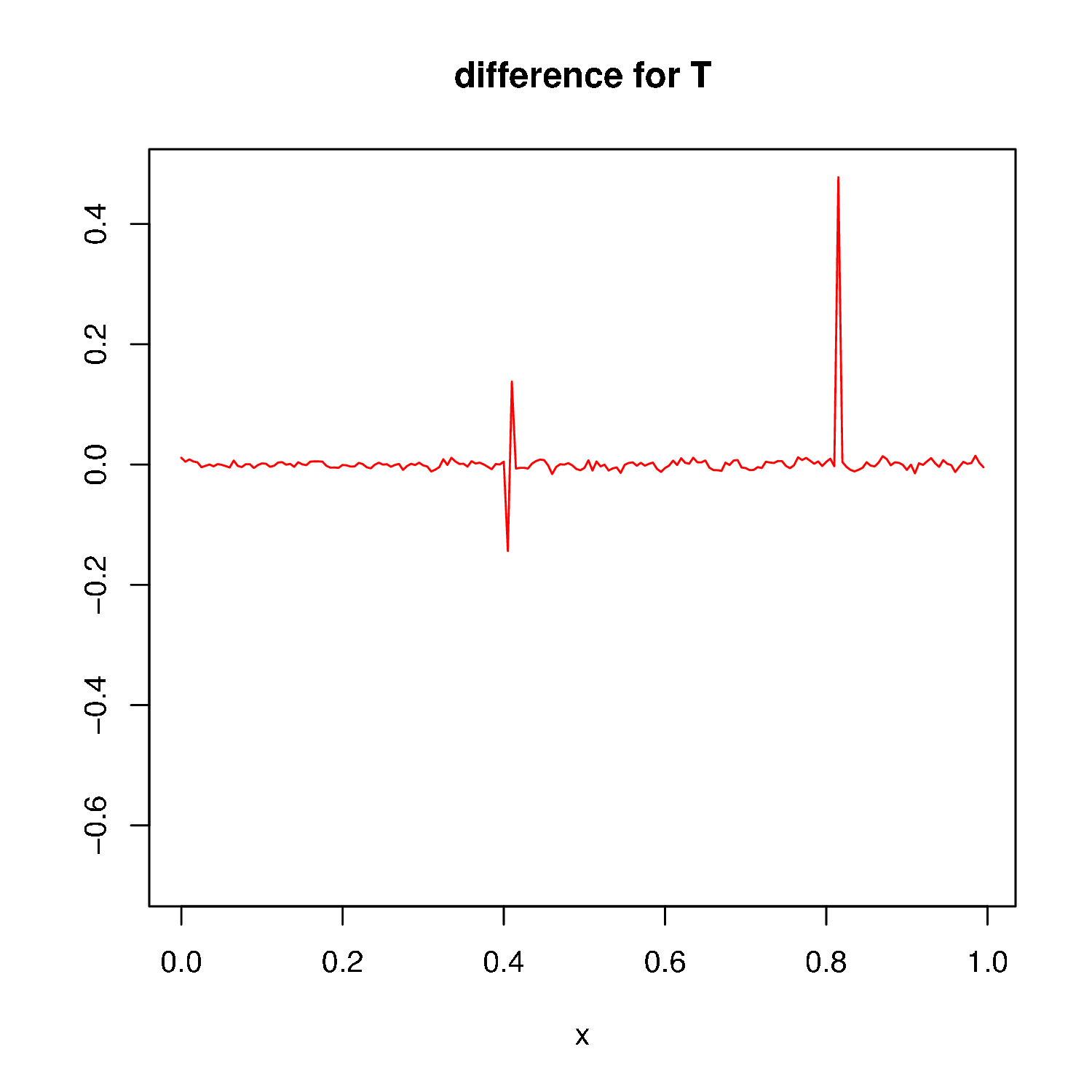}
\end{center}
\subsubsection{Infinite number of pieces}
If the number of pieces is infinite, their are many different settings leading to existence and uniqueness of an absolutely continuous invariant probability measure which is exponentially mixing. Let us cite \cite{Br,LSV}. We will consider the conditions of Liverani-Saussol-Vaienti \cite{LSV} in a restricted case (our potential is $|T'|^{-1}$, they consider more general potential). These conditions are the following.
\begin{cond}\label{infinite}
\begin{enumerate} 
\item  There is a countable partition $(I_j)_{j\in\N}$ of $I$ into intervals. On each $I_j$, the map $T$ is monotonic and $C^2$, it is continuous on each $\overline{I_j}$. Denote by $\PP$ the partition  $(I_j)_{j\in\N}$ and by $\PP_k$ the partition of monotonicity of $T^k$.
\item $\displaystyle \frac1{|T'|}\in BV$ and $\displaystyle \sum_{P\in\PP}\sup  \frac1{|T'|}<\infty$.
\item  The partition $\PP$ is generating.
\item The map $T$ is covering, which means, in the infinite case, that for any $n\in\N$, for any $P\in\PP_n$, $I$  may covered by a finite number  of subintervals of $T^N P$~:
$$\forall n\in \N \/, \ \forall P\in\PP_n \/, \ \exists N \/, \ \exists \ \mbox{finite} \ {\mathcal J}\subset \PP_N\bigvee \{P\} \ / \ \bigcup_{Q\in{\mathcal J}}T^NQ =X\/.$$
\end{enumerate}
\end{cond}
\begin{rem}
The condition  $\displaystyle \frac1{|T'|}\in BV$ is satisfied provided that $T$ has the bounded distortion property~:
$$\sup_{P\in \PP}\sup_{x\in P} \frac{|T''(x)|}{|T'(x)|^2} <\infty \/$$
and that $\displaystyle \sum_{P\in\PP}\sup  \frac1{|T'|}<\infty$.
\end{rem}
The following result has been proved in a general setting in \cite{LSV} and in \cite{Br} under the condition of bounded distortion.
\begin{theo}\label{vitinfinite}
Under Assumptions \ref{infinite}, the map $T$ admits a unique absolutely invariant probability measure, its density $f$ belongs to $BV$ and $\inf f >0$. Moreover, if $\mu = f m$ is this invariant probability measure, then (\ref{mixing}) is satisfied with $\CC = BV$ and $\Phi(n)= \gamma^n$, $0<\gamma<1$~: there exists $0<\gamma<1$, $C>0$, such that for any $\psi\in BV$ and $\varphi \in L^1(\mu)$, for any $n\in\N$, 
$$\left|\int\limits_I \psi \cdot \varphi\circ T^n \/ d\mu - \int\limits_I \psi \/ d\mu \int\limits_I \varphi \/ d\mu \right| \leq C \/ \gamma^n\/ \Vert \varphi \Vert_1 \/ \Vert \psi \Vert_{BV}\/.$$
\end{theo}
\begin{lem}
The map $T$ satisfying Assumptions \ref{infinite} is regular.
\end{lem}
\begin{proof}[Idea of the proof]
This is a simple consequence of the fact that $T$ is piecewise $C^2$ (thus piecewise $C^1$) and that $\displaystyle \sum_{P\in\PP}\sup  \frac1{|T'|}<\infty$.
\end{proof}
There are several natural examples of dynamical systems satisfying Assumptions \ref{infinite}. Let us cite the Gauss map~:
$$T(x) = \frac1x -\left\lfloor \frac1x \right\rfloor $$
which appears in the continuous fractions decomposition (\cite{Br}) and in analysis a gcd algorithms (\cite{V}). An natural extension of these systems are ``Japanese systems'' or $\alpha$-Gauss maps~: for $0<\alpha\leq 1$,
$$T_\alpha(x) =\left|\frac1x\right|-\left\lfloor \left|\frac1x\right|+1-\alpha \right\rfloor\/,$$ 
$T_\alpha$ maps the interval $[\alpha-1\/, \alpha]$ into itself. See \cite{BDV} for a description of these systems as well as an application to analysis of generalized Euclidian algorithms. The maps $T_\alpha$ satisfy Asumption \ref{infinite} for $0<\alpha \leq 1$.
\begin{coro}\label{infinitepieces}
Let $T$ satisfies Assumption \ref{infinite}, let $K$ be a Kernel belonging to $BV$, let $\hat{f}_N$ and $\hat{T}_N$ be the estimators of $f$ and $T$.  For all $0<\alpha<1$, here exists $M>0$, $L>0$, $R>0$,   such that outside a set of measure less than $R h^{\alpha}$,   for all $t\in\R^{+}$,   for all $u\geq \mbox{Ct} \/h$, 
\begin{align}
\P(|\hat{f}_n(x)-f(x)|>t-u^\alpha)\leq 2e^{\frac1e} \exp[-t^2 M nh^2] \\
\P(|\hat{T}_n(x)-T(x)|>t-u^\alpha) \leq 2e^{\frac1e} \exp[-t^2 L  nh^2] \/.
\end{align}
As a consequence,  provided $h=h_n$ goes to zero  and $n h_n^2 =O(n^\eps)$, $\eps>0$, we obtain the following convergences~:
\begin{itemize}
\item for $m$ almost all $x\in I$, $\hat{f}_n(x)$ converges to $f(x)$ and $\hat{T}_n(x)$ converges to $T(x)$ almost surely and in $L^p$  for any $1\leq p$,
\item $\displaystyle\E\left(\int_I |\hat{f}_n(x) -f(x)| \/dx \right)$ and $\displaystyle\E\left(\int_I |\hat{T}_n(x)-T(x)|\right)$ go to zero.
\item for almost all $x\in I$, for $a <\frac12$, $|\hat{g}_n(x)-g(x)| =O_\P(n^{-a})$ where $\hat{g}_n$ is either $\hat{f}_n$ or $\hat{T}_n$ and $g$ is either $f$ or $T$.
\end{itemize}
\end{coro}
The following graphics are for the Gauss map $T(x)=\frac1x \ (\mbox{mod} \ 1)$ for which it is known that the unique invariant probability density is $f(x) = \frac1{\log 2} \frac1{1+x}$. We have performed simulations with $n=50 000$, $h=0.009$ and the noise $\eps_i\leadsto {\mathcal U}[-0.2\/,0.2]$. Over $800$ points, we get $EMEf=  0.05046933$, $AMET=0.05141938$, if we restrict to points in $[0.2\/,1]$ we get $AMET=0.01439787$.
\begin{center}
\includegraphics[scale=0.43]{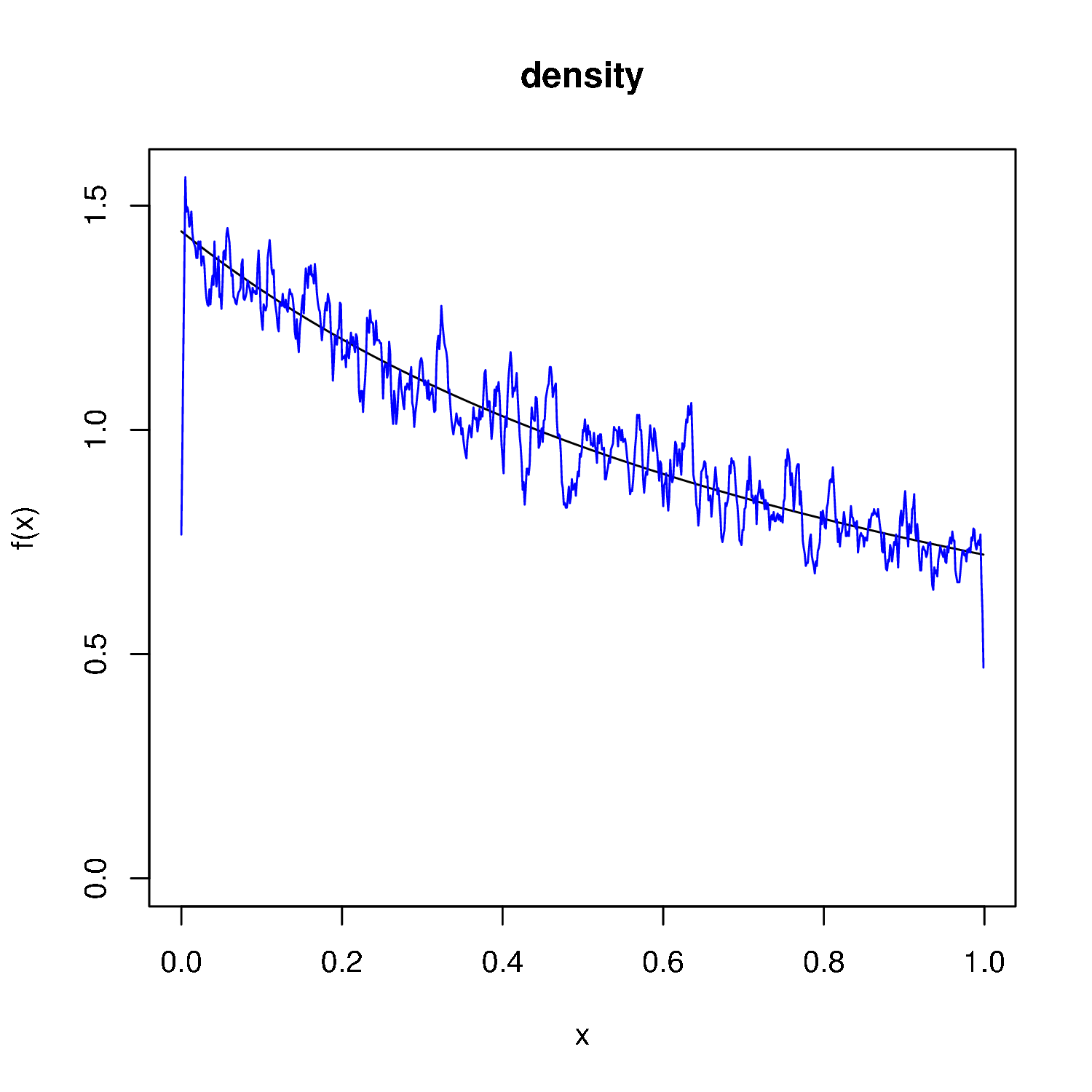}
\includegraphics[scale=0.43]{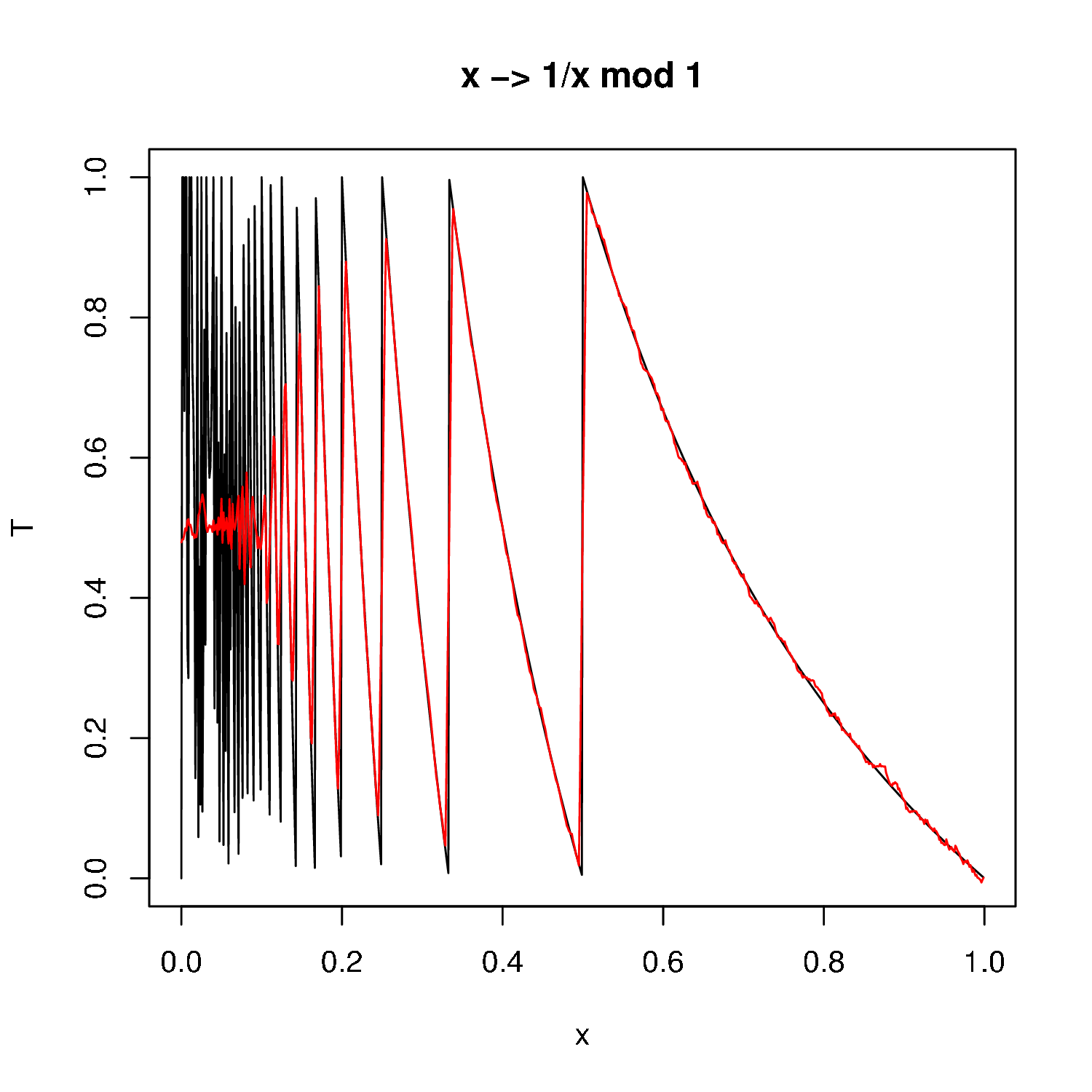}
\includegraphics[scale=0.43]{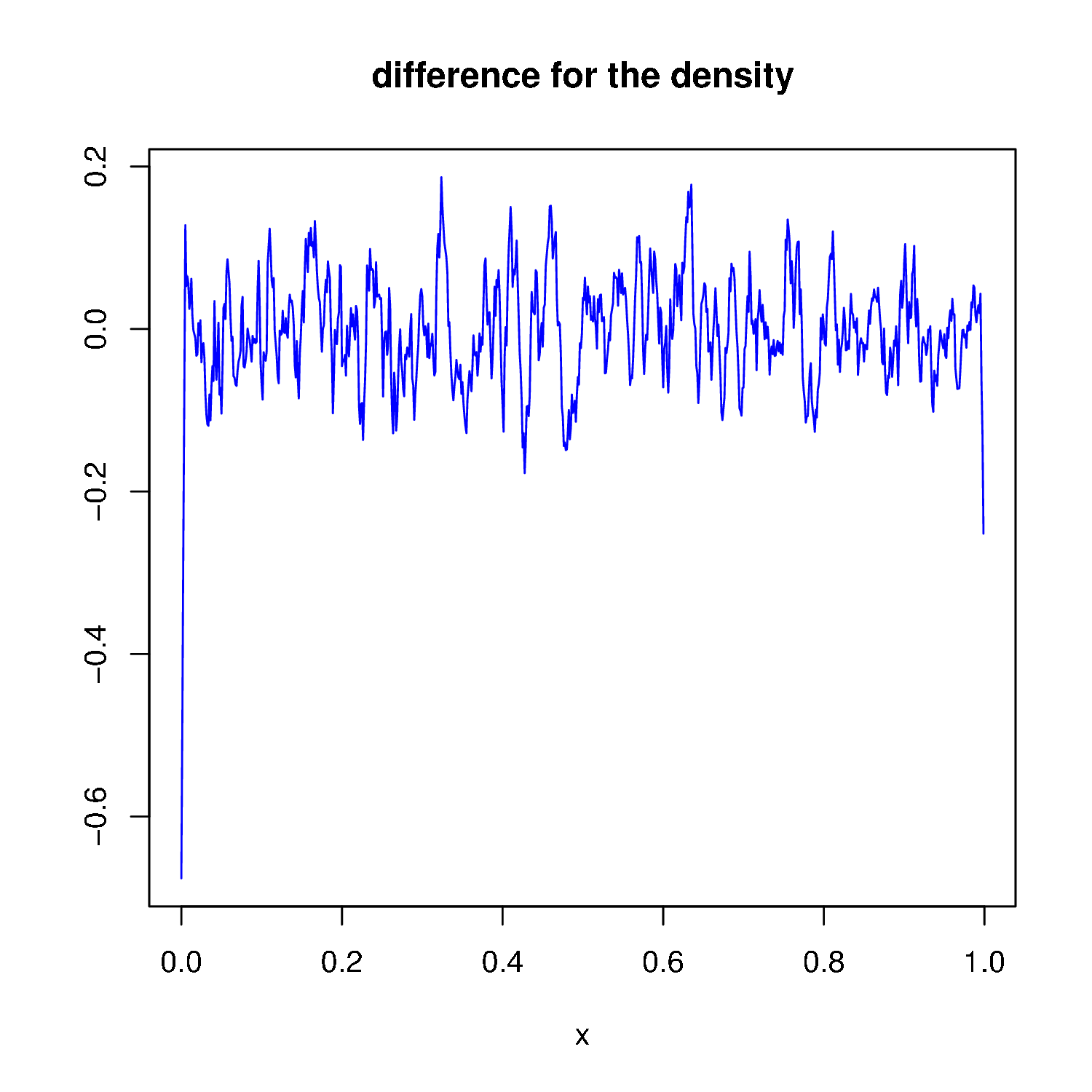}
\includegraphics[scale=0.43]{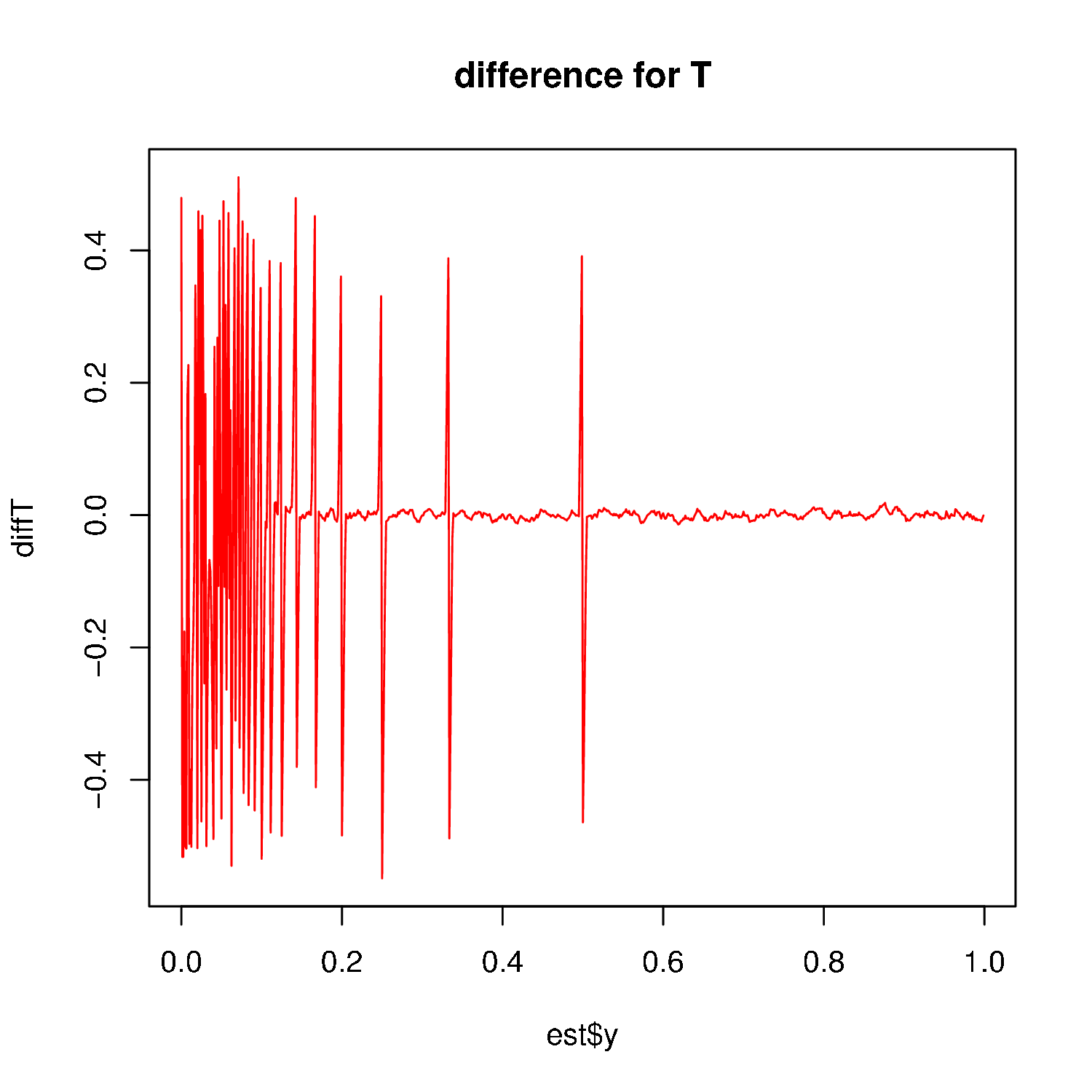}
\end{center}
\subsection{In dimension one~: unimodal maps}
Let $I=[-1, 1]$ and $f: I \to I$ be a $C^2$ unimodal map (i.e., $f$ is increasing on $[-1,0]$, decreasing on   $[0,1]$) satisfying $f''(0) \ne 0$, and,   \smallskip   
     
\noindent {\bf (H1)} There are   $0 < \alpha < 1$, $K > 1$, and    $\tilde \lambda  \leq \lambda \le 4$ with  $e^{ 2\alpha}< \tilde \lambda $, and $\sup_I|f'|\le \tilde \lambda^{K}< 8$ so that \\  
{(i)}  $|(f^n)'(f(0))|\ge \tilde \lambda^n$ for all  $n \in \N$ and   $\lambda = \lim_{n \to \infty} |(f^n)'(f(0))|^{1/n} $.  \\
{(ii)} $|f^n(0)| \ge  e^{-\alpha n}$, for all $n \ge 1$.  
  

\noindent {\bf (H2)} For each small enough  $\delta > 0$,  there is $M=M (\delta) \in \N_+$ for which   \\ 
{(i)} If $x , \ldots, f^{M-1} (x)    \notin (-\delta, \delta)$ then $|(f^M)'(x)| \ge \tilde \lambda ^M$;      
{(ii)} For each $n$, if $x, \ldots, f^{n-1}(x) \notin   (-\delta, \delta)$ and $f^n(x) \in (-\delta, \delta)$,  then $|(f^n)'(x)| \ge \tilde \lambda^n$.   \\
\noindent {\bf (H3)}   $f$ is topologically mixing on $[f^2(0), f(0)]$, that is for any two open sets $U$, $V\subset I$, there exists $N\in\N$ such that $\forall n\geq N$, $T^{-n}U\cap V\neq \emptyset$.   
   
\medskip   
Examples of unimodal maps satisfying (H1)--(H3) are quadratic maps $1-a\cdot x^2$ for a positive  measure set of parameters $a$.  (See e.g. \cite{BeCa}).  \\
The following theorem is obtained in two steps. First, it is proven that unimodal maps satisfying (H1)--(H3) are conjugated to a kind of hyperbolic markov maps called Young towers (see \cite{Y1,KN}). Then the estimation on the speed of mixing is obtained on the tower (\cite{BuM1}). Let us emphasize the fact that the kind a mixing we need (namely \ref{mixing}) and especially the fact that  $\Vert \varphi \Vert_1$ appears, is not obtained in \cite{Y1,KN} (a form with the $\Vert \varphi \Vert_\infty$ was obtained there). 
\begin{theo}
Let $T$ be a unimodal map satisfying (H1)--(H3).  The map $T$ admits a unique absolutely invariant probability measure with density $f$ satisfying $\displaystyle \inf_{[f^2(0)\/,f(0)]} f >0$. Moreover, if $\mu = f m$ is this invariant probability measure, then (\ref{mixing}) is satisfied with $\CC$ the space  $\Lip$ of Lipschitz functions on $I$ and $\Phi(n)= \gamma^n$, $0<\gamma<1$~: there exists $0<\gamma<1$, $C>0$, such that for any $\psi\in \Lip$ and $\varphi \in L^1(\mu)$, for any $n\in\N$, 
$$\left|\int\limits_I \psi \cdot \varphi\circ T^n \/ d\mu - \int\limits_I \psi \/ d\mu \int\limits_I \varphi \/ d\mu \right| \leq C \/ \gamma^n\/ \Vert \varphi \Vert_1 \/ \Vert \psi \Vert_{\Lip}\/.$$
\end{theo}
This is clear that $T$ is $1$-regular~: for well chosen $\eps_n$ and $u_n$, the sets  $A_N$ and $B_N$ are empty because $T$ is $C^2$. Nevertheless, it is known that the invariant density is very irregular (see \cite{Y2,KN}). With a more intricate study, we  should probably prove that $f$ is also regular. Here, we restrict ourselves to the estimation of $T$. \\
Since the invariant measure has its support in $S=[f^2(0)\/,f(0)]$, our estimates are valid only for $x \in S$.
\begin{coro}\label{unimodales}
Let $T$  be a unimodal map satisfying (H1)--(H3), let $K$ be a Kernel belonging to $\Lip$, let $\hat{T}_N$ be the estimator of  $T$. There exists $M>0$, $L>0$, $R>0$,   such that outside a subset of $S$, of measure less than $R h$,   for all $t\in\R^{+}$,   for all $u\geq \mbox{Ct} \/h$, 
$$\P(|\hat{T}_n(x)-T(x)|>t-u) \leq 2e^{\frac1e} \exp[-t^2 L h^2  n] \/.$$
As a consequence,  provided $h=h_n$ goes to zero  and $n h_n^2 =O(n^\eps)$, $\eps>0$, we obtain the following convergences~:
\begin{itemize}
\item for $m$ almost all $x\in S$,  $\hat{T}_n(x)$ converges to $T(x)$ almost surely and in $L^p$  for any $1\leq p$,
\item $\displaystyle\E\left(\int_S |\hat{T}_n(x)-T(x)|\right)$ go to zero.
\item for almost all $x\in S$, for $a <\frac12$, $|\hat{T}_n(x)-T(x)| =O_\P(n^{-a})$.
\end{itemize}
\end{coro}
The following graphics are for the map $T(x)=3.8*x*(1-x)$. We have performed simulations with $n=50 000$, $h=0.01$ and the noise $\eps_i\leadsto {\mathcal U}[-0.2\/,0.2]$. Over $154$ points in $S$, we get  $AMET=0.004114143$.
\begin{center}
\includegraphics[scale=0.43]{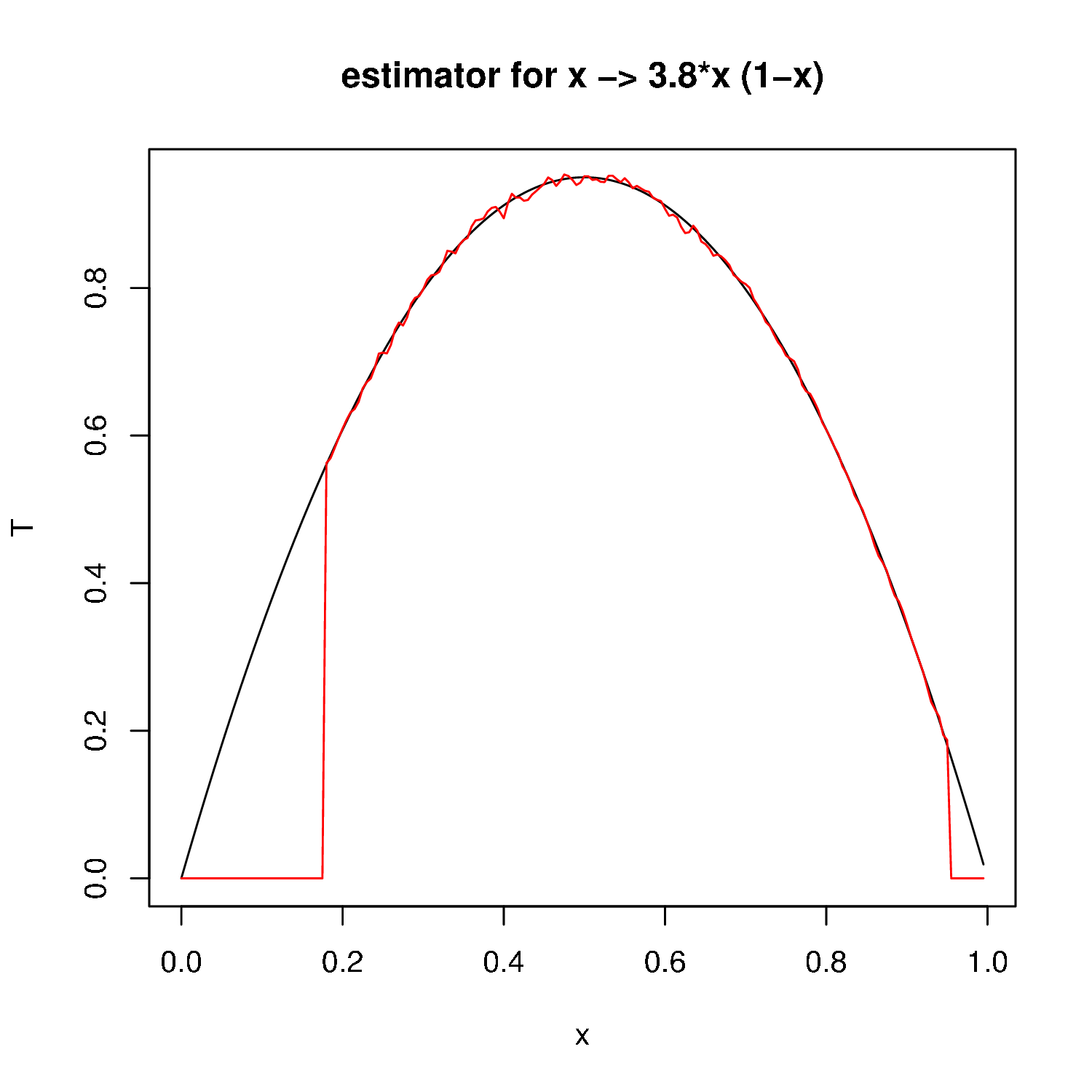}
\includegraphics[scale=0.43]{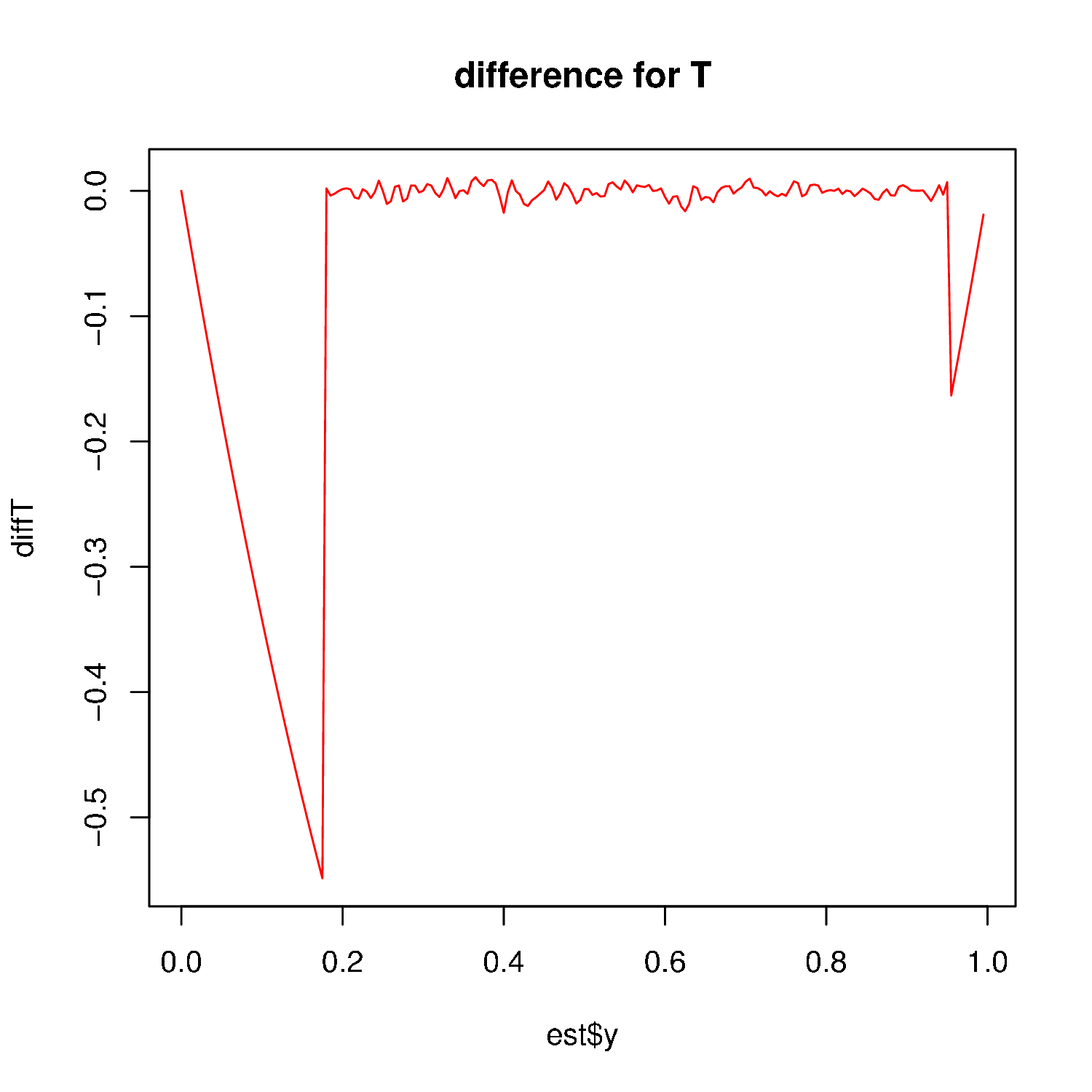}
\end{center}
\subsection{Piecewise expanding maps in higher dimension} These are generalisations of Lasota-Yorke maps in higher dimension. These maps have been studied in \cite{Bu2,BPS,Bu0,Buz,Bu3,Cow,GB,Sau} from various point of view. The control of the speed of mixing may be found in \cite{BuM1,BuM2}, the strategy is as for unimodal maps~: the map is conjugated to a ``Young tower''.\\
The setting is  the following. $(X,\ZZ,T)$ will be a {\em  piecewise invertible map}, i.e.:
\smallskip
\begin{itemize}
  \item $X=\overline{\bigcup_{Z\in\ZZ} Z}$ is a locally connected compact subset of $\R^d$. 
  \item $\ZZ$ is a finite collection of pairwise disjoint, bounded and open
subsets of $X$, each with a non-empty boundary. Let $Y=\bigcup_{Z\in\ZZ} Z$.
  \item $T:Y\to X$ is a map such that each restriction $T|Z$,   $Z\in\ZZ$, coincides with the restriction of a 
homeomorphism $T_Z:U\to V$ with $U,V$ open sets such that $U\supset\bar{Z}$, $V\supset\overline{T(Z)}$.
\end{itemize}
$T$ will be assumed to be {\em non-contracting}, i.e., such that for all $x,y$ in the same element $Z\in\ZZ$, $d(Tx,Ty)\geq d(x,y)$. Also $\ZZ$ will be assumed to be {\em generating}, i.e., $\lim_{n\to\infty} \mbox{diam}(\ZZ^n)=0$ where $\ZZ^n$ denotes the set of $n$-cylinders, i.e., the non-empty sets of the form:
 $$
     [A_0\dots A_{n-1}] := A_0\cap\dots\cap T^{-n+1}A_{n-1}
 $$
for $A_0,\dots,A_{n-1}\in\ZZ$. 
\medskip
Finally the {\em boundary of the partition}, $\partial\ZZ=\bigcup_{Z\in\ZZ} \partial Z$, will play an important role in our analysis. In particular, we shall assume ``small boundary pressure'' (see below), a fundamental condition which already appeared in \cite{Bu2,Bu3,BPS}.
\medskip 
To formulate the crucial ``small boundary pressure'' condition, we need first some definitions.\\
The {\em topological pressure} \cite{DKS} of a subset $S$ of $X$ is:
$$
   P(S,T) = \limsup_{n\to\infty} \frac1n \log
      \sum_{{A\in\ZZ^n}\atop{ \bar{A}\cap S\ne\emptyset}}
     g^{(n)}(A)
$$
where $g^{(n)}(A)=\sup_{x\in A} g(x)g(Tx)\dots g(T^{n-1}x)$, $g= |\det T'|^{-1}$.\\
The {\em small boundary pressure} condition is:
$$   P(\partial\ZZ,T) < P(X,T). $$
\medskip
This inequality is satisfied in many cases. In particular, if $T$ is expanding and $X$ is a Riemannian manifold and the weight is $|\det T'(x)|^{-1}$ or close to it, then it is satisfied: (i) in dimension $1$, in all cases;   (ii) in dimension $2$, if $T$ is piecewise real analytic \cite{Bu3,Tsu1}; (iii) in arbitrary dimension, for all piecewise affine $T$ \cite{Tsu2} or for generic $T$ \cite{Bu1,Cow}. \\
A basic example is given by the {\em multidimensional $\beta$-transformations} \cite{Bu0}, i.e., maps $T:[0,1]^d\to[0,1]^d$, $T(x)=B(x)\; {\rm mod}\; \Z^d$ with $B$ an expanding affine map on $\R^d$.
Let us summarize our hypothesis on $T$.
\begin{cond}\label{dimsup}
Let $(X,\ZZ,T,g)$ be a weighted piecewise invertible dynamical system.
Assume that:
\begin{itemize}
\item $T$ is expanding, i.e., there is some $\lambda>1$ such that for all $x,y$ in the same element of $\ZZ$, $d(Tx,Ty)\geq \lambda\cdot d(x,y)$;
\item $g= |\det T'|^{-1}$ is H{\"o}lder continuous with exponent $\gamma$ and is
positively lower bounded.
\item the boundary pressure is small: $P(\partial\ZZ,T)<P(X,T)$.
\item $T$ is topologically mixing.
\end{itemize} 
\end{cond}
Let $K(f)=\max_{Z\in\ZZ} \sup_{x\ne y\in Z}
\frac{|f(x)-f(y)|}{d(x,y)^\gamma}$ where $\gamma$ is some H{\"o}lder
exponent of $g$.
\begin{theo}\cite{BuM1}
Let $T$ satisfy Assumption \ref{dimsup}. Then,  $T$ admits a unique invariant measure $\mu$, absolutely continuous w.r.t. the Lebesgue meausre $m$. This measure is exponentially mixing~:
$$
   \left| \int_X \varphi\circ T^{n} \cdot \psi \, d\mu 
       - \int_X \varphi \, d\mu\int_X \psi \, d\mu \right|
      \leq C \cdot \|\varphi\|_{C^\gamma(X)} \cdot \|\psi\|_{L^1}
   \; \kappa ^n
\/.
$$
with constants $C<\infty$ and $\kappa<1$ depending only on
$(X,\ZZ,T,g)$, for any measurable functions
$\varphi,\psi:X\to\R$ such that $\psi$ is bounded and $\varphi$ is
$\gamma$-H{\"o}lder continuous. 
\end{theo}
This is clear that $T$ is $1$-regular because it is piecewise $C^1$. Nevertheless, it is known that the invariant density has discontinuities on $\partial [T^n(\ZZ^n)]$ . With a more intricate study, we  should probably prove that $f$ is also regular. Here, we restrict ourselves to the estimate of $T$.
\begin{coro}\label{resdimsup}
Let $T$ satisfy Asumption \ref{dimsup}, let $K$ be a $\gamma$-H\"older Kernel, let $\hat{T}_N$ be the estimator of  $T$. There exists $M>0$, $L>0$, $R>0$,   such that outside a set of measure less than $R h$,   for all $t\in\R^{+}$,   for all $u\geq \mbox{Ct} \/h$, 
$$\P(|\hat{T}_n(x)-T(x)|>t-u) \leq 2e^{\frac1e} \exp[-t^2 L h^{\gamma+2}  n] \/.$$
As a consequence,  provided $h=h_n$ goes to zero  and $n h_n^{\gamma+2} =O(n^\eps)$, $\eps>0$, we obtain the following convergences~:
\begin{itemize}
\item for $m$ almost all $x\in X$,  $\hat{T}_n(x)$ converges to $T(x)$ almost surely and in $L^p$  for any $1\leq p$,
\item $\displaystyle\E\left(\int_X |\hat{T}_n(x)-T(x)|\right)$ go to zero.
\item for almost all $x\in X$, for $a <\frac12$, $|\hat{T}_n(x)-T(x)| =O_\P(n^{-a})$.
\end{itemize}
\end{coro}
We have performed some simulations for $T(x) = Bx \  \mbox{mod}  \ \Z^2$, with $B$ the matrix 
$$\left( \begin{matrix}
2.5 & 3.4 \\
4.6 & 3.2 
\end{matrix}\right)$$
Below are the histograms for $diff_x$ (resp. $diff_y$), the difference beetween the $x$ (resp. $y$) coordinate of $T$ and the $x$ (resp. $y$) coordinate of $\hat{T}_n$, over a grid of $100$ times $100$ points in $[0\/,1]^2$. The Kernel is $K = \frac14\Id_{[-1\/,1]\times [-1\/,1]}$, there is no noise, $\eps_i =0$, $n=66668$ and $h=0.004$. The AME for the coordinate $x$ is $0.01882885$, for the coordinate $y$, the AME is $0.06723186$.
\begin{center}
\includegraphics[scale=0.43]{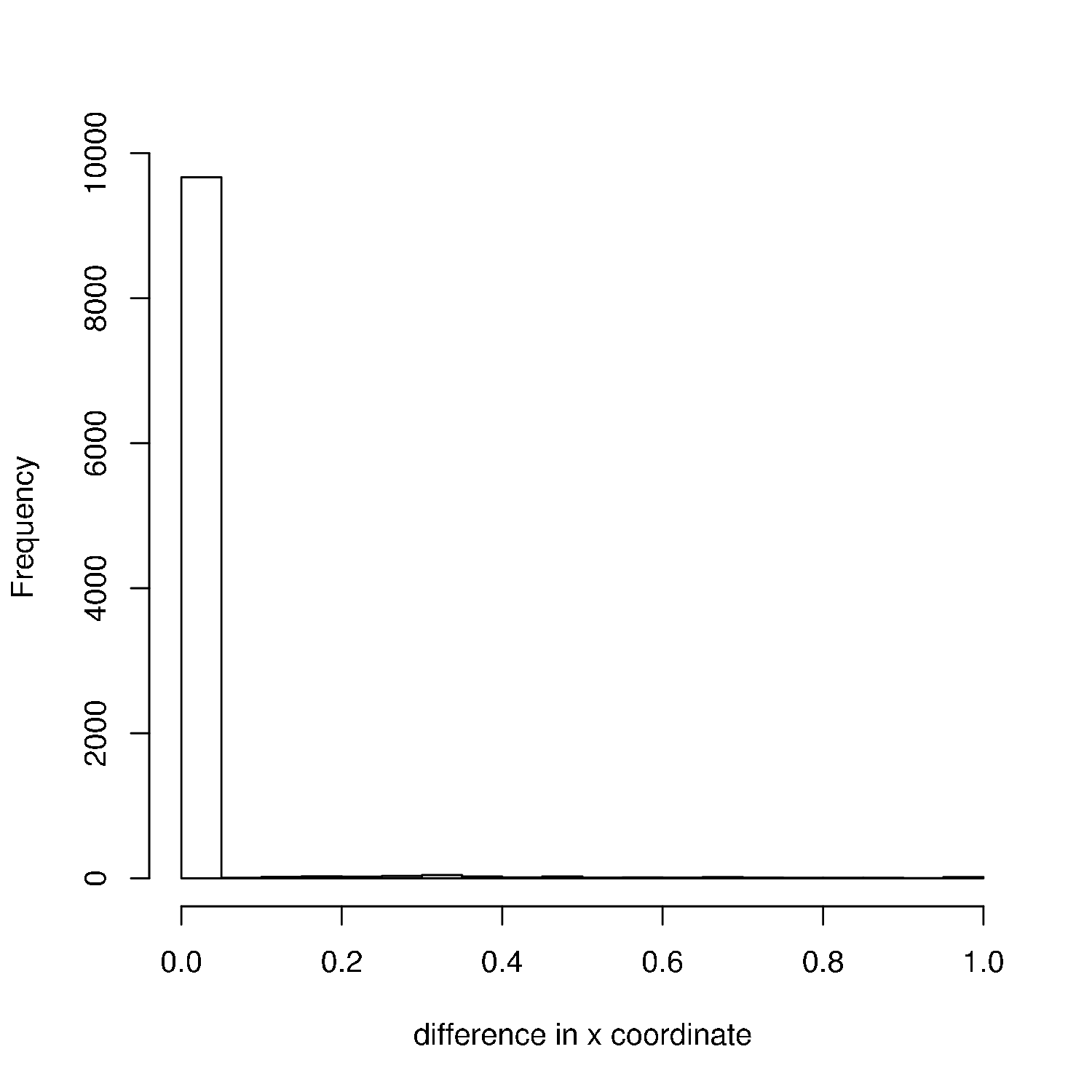}
\includegraphics[scale=0.43]{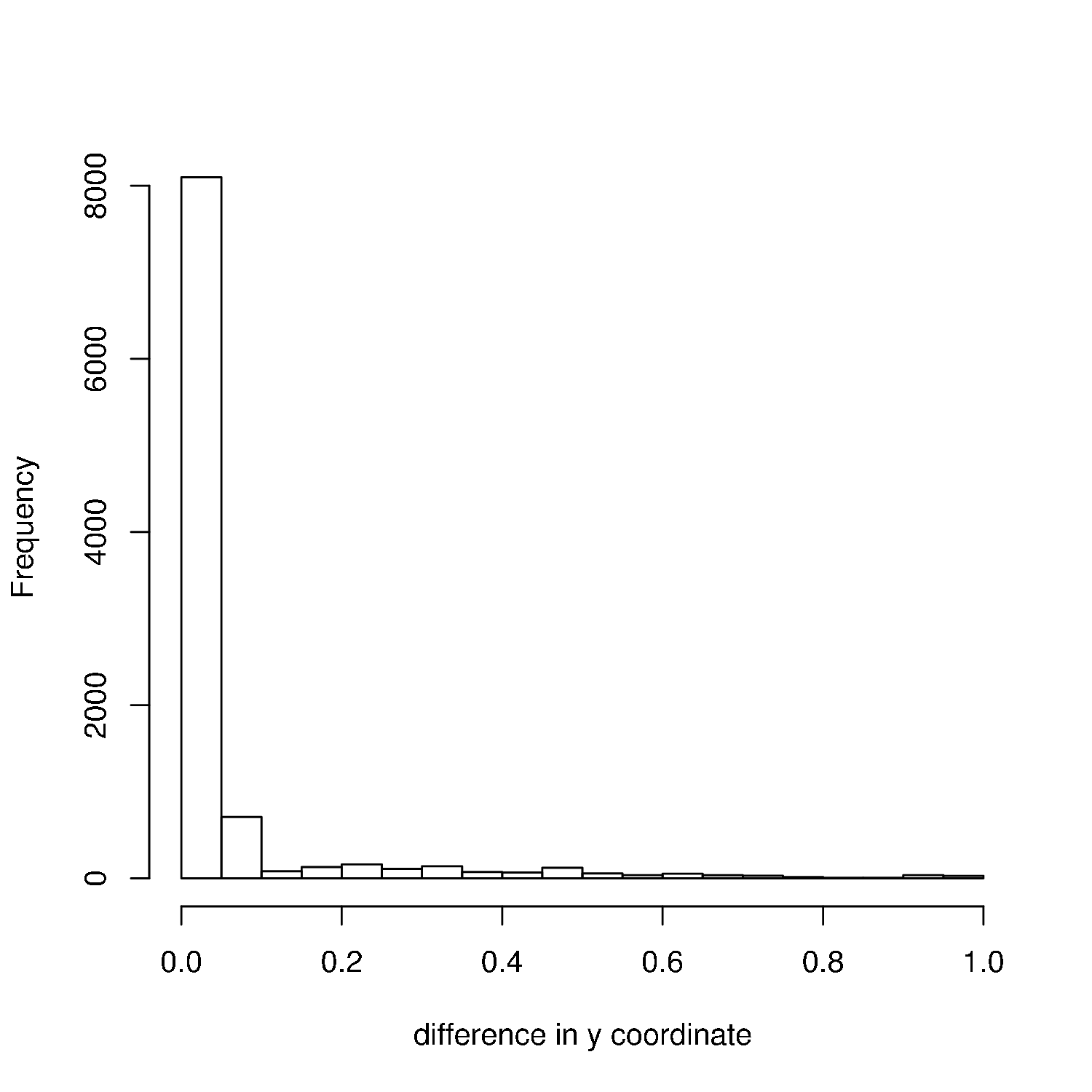}
\end{center}
\begin{rem}[Anosov maps] 
Our technics should also apply to estimate the invariant density and the application $T$ for Anosov maps for which there exists an invariant measure absolutely continuous with respect to the Lebesgue measure. A more intricate study could also lead results of the same kind for Axiom A diffeomorphisms.
\end{rem}

\end{document}